\documentclass[10pt,a4paper]{scrartcl}

\usepackage{a4wide}
\usepackage{latexsym}
\usepackage{amssymb}
\usepackage{mathrsfs}
\usepackage{graphicx}
\usepackage[utf8]{inputenc}
\usepackage{amssymb}
\usepackage{amsmath}
\usepackage{amsfonts}
\usepackage{amsthm}
\usepackage{amsxtra}
\usepackage{url}
\usepackage{mathptmx}
\usepackage{helvet}
\usepackage{stmaryrd}
\usepackage{dsfont}
\usepackage{setspace}
\usepackage{hyperref}
\usepackage{enumerate}
\usepackage{MathDef}
\usepackage{verbatim}
\usepackage{color}
\usepackage{cleveref}
\usepackage[comma]{natbib}

\numberwithin{equation}{section}
\allowdisplaybreaks

\newtheorem{Theorem}{Theorem}[section]
\newtheorem{Corollary}{Corollary}[section]
\newtheorem{Proposition}{Proposition}[section]
\newtheorem{Lemma}{Lemma}[section]
\newtheorem{Definition}{Definition}[section]

\theoremstyle{definition}
\newtheorem{Remark}{Remark}[section]

\usepackage{titlesec}
\titleformat{\paragraph}{\normalfont\bfseries}{\theparagraph}{}{}
\titlespacing*{\paragraph}{0pt}{3.25ex plus 1ex minus .2ex}{0.1em}

\date{}

\title{\LARGE Dependence Estimation for  High Frequency\\[2mm] Sampled Multivariate CARMA Models}

\author{Vicky Fasen }

\date{}

\begin{document}
%
\maketitle
%
\begin{abstract}
The paper considers high frequency sampled multivariate continuous-time ARMA \linebreak (MCARMA)
models, and derives the
asymptotic behavior of the sample autocovariance function to a normal random matrix.
Moreover, we obtain the asymptotic behavior of the cross-covariances between different
components of the model. We will see that the limit distribution of the sample autocovariance function
has a similar structure in the continuous-time and in the
discrete-time model.
As special case we consider a CARMA (one-dimensional MCARMA)
process. For a CARMA process we prove Bartlett's formula for the
sample autocorrelation  function. Bartlett's formula has the same form in both models,
only the sums in the discrete-time model are exchanged by integrals in the continuous-time model.
Finally, we present limit results for multivariate MA processes as well
which are not known in this generality in the multivariate setting yet.

\end{abstract}

\vspace{0.2cm}

\noindent
\begin{tabbing}
\\ \> Secondary: 60F05,\,60G10 
\end{tabbing}

\vspace{0.2cm}\noindent\emph{Keywords:}   asymptotic normality,  Bartlett's formula, CARMA process, consistency, correlation function, covariance function,
high-frequency data, L\'evy process, limit theorems, MA process, multivariate models, VARMA process.

\section{Introduction}

The paper considers multivariate ARMA (autoregressive moving average) models
in continuous time and their dependence estimation.
Multivariate time series have the advantage that they are able to model dependence
between different time series in the most generality.
 A classical dependence measure for a multivariate stationary
process $(\Ybf_t)_{t\in\R}$ is the autocovariance function. The autocovariance function
is defined as $\Gamma_\Ybf(h)=\E((\Ybf_0-\E(\Ybf_0))(\Ybf_h-\E(\Ybf_h))^T)$, $h\in\R$.
An estimator for the autocovariance function is the sample autocovariance function.

One of the most known multivariate time series models is the  {\em VARMA$(p,q)$} (vector autoregressive moving average) process $(p,q\in\N_0)$
defined to be the stationary solution to a $d$-dimensional difference equation of the form
\begin{eqnarray} \label{VARMA}
    \mathbf{P}(B)\Ybf_k=\mathbf{Q}(B)\xi_k,
\end{eqnarray}
where $B$ is the backshift operator satisfying $\Bbf\Ybf_k=\Ybf_{k-1}$, $(\xi_k)_{k\in\Z}$
is a sequence of independent and identically distributed (iid) random vectors in $\R^m$,
\beam \label{Pol P}
    \mathbf{P}(z):=I_{d}z^p+\Pbf_1z^{p-1}+\ldots+\Pbf_{p-1}z+\Pbf_p
\eeam
with $\Pbf_1,\ldots, \Pbf_p\in \R^{d\times d}$ is the autoregressive polynomial and
\beam \label{Pol Q}
    \mathbf{Q}(z):=\Qbf_0z^q+\Qbf_1z^{q-1}+\ldots+\Qbf_{q-1}z+\Qbf_q
\eeam
with $\Qbf_0,\ldots,\Qbf_q\in \R^{d\times m}$
is the moving average polynomial. In this article we always assume that \linebreak $\E\|\xi_1\|^2<\infty$.
If $\det(\Pbf(z))\not=0$ for all $z\in\C$ with $|z|\leq 1$, then \eqref{VARMA} has exactly
one solution which has the moving average representation
$    \Ybf_k=\sum_{j=0}^\infty\Cbf_j\xi_{k-j}$,
where the $\Cbf_j$ are uniquely determined by $\Cbf(z)=\sum_{j=0}^\infty\Cbf_j z^j=\Pbf(z)^{-1}\Qbf(z)$
for $|z|\leq 1$ (cf. \cite{Brockwelletal1991}, Theorem 11.3.1). The interest in the asymptotic properties
of the sample autocovariance and the autocorrelation function
of ARMA (one-dimensional VARMA) processes has a long history starting with
\cite{Bartlett1955} and continuing with \cite{Anderson:Walker:1964,Hannan76,Cadena94} to name a few (cf. \cite{Brockwelletal1991}).
However, for VARMA processes the multivariate nature of the covariance matrix $\Gamma(h)$ is a challenge, and hence, for a very long time people looked only at special cases like the asymptotic
behavior of cross-covariances of bivariate Gaussian MA processes or independent MA processes (cf. \cite{Brockwelletal1991,Fuller:book}).
Quite recently  \cite{Su:Lund} developed the asymptotic behavior of multivariate MA processes in a more general setup.

{\em Multivariate continuous-time ARMA} (MCARMA) processes $\Ybf=(\Ybf_t)_{t\in\R}$
are the continuous-time
versions of VARMA processes.
The driving force of a MCARMA process is a {\em L\'{e}vy process} $(\Lbf_t)_{t\in\bbr}$. A L\'{e}vy process $(\Lbf_t)_{t\geq 0}$  is defined
to satisfy $\Lbf_0=0$ a.s., $(\Lbf_t)_{t\geq 0}$ has independent and stationary increments and
the paths of $(\Lbf_t)_{t\geq 0}$ are stochastically continuous. An extension of a L\'{e}vy process $(\Lbf_t)_{t\geq 0}$ from the
positive to the whole real line
is given by $\Lbf_t := \Lbf_t\mathds{1}_{\{t\geq 0\}} - \wt \Lbf_{-t-}\mathds{1}_{\{t<0\}}$ for $t\in\bbr$, where
$ (\wt \Lbf_t)_{t\geq 0}$  is an independent copy of $(\Lbf_t)_{t\geq 0}$.
Prominent examples are Brownian motions, compound Poisson processes and stable L\'{e}vy processes.
L\'{e}vy processes  are
characterized by their {\em L\'{e}vy-Khintchine representation}.
An $\R^m$-valued L\'{e}vy process $(\Lbf_t)_{t\geq 0}$ has the
 L\'{e}vy-Khintchine representation
$\E(\e^{i \Theta^T\Lbf_t})=\exp(-t\Psi(\Theta))$
for $\Theta\in\mathbb{R}^m$
and
 \beao
        \Psi(\Theta)=
        -i \gabf_{\Lbf}^T\Theta
        +\frac{1}{2} \Theta^T\Sibf_{\Lbf}\,\Theta +\int_{\R^m}
 \left(1-\e^{i \xbf^T\Theta}+i \xbf^T\Theta\1_{\{\|\xbf\|^2\leq 1\}} \right)\nu_{\Lbf}(\dd\xbf)
    \eeao
with  $\gabf_{\Lbf}\in\mathbb{R}^m$, $\Sibf_{\Lbf}$ a positive semi-definite
matrix in $\mathbb{R}^{m\times m}$ and
 $\nu_{\Lbf}$  a measure on $(\mathbb{R}^m,\mathcal{B}(\mathbb{R}^m))$,
 called  {\em L\'{e}vy measure}, which satisfies
$\int_{\R^m}\min\{\|\xbf\|^2,1\}\,\nu_{\Lbf}(\dd \xbf)<\infty$
and  $\nu_{\Lbf}(\{\0_m\})=0$. The triplet
$(\gabf_{\Lbf},\Sibf_{\Lbf},\nu_{\Lbf})$ is called
the {\em characteristic triplet}, because
it characterizes completely the distribution of the L\'{e}vy process. For more details on L\'{e}vy
processes we refer to the excellent monograph of \cite{Sato1999}.

Let $\Lbf=(\Lbf_t)_{t\in\R}$ be a two-sided $\R^m$-valued L\'{e}vy process and let $p>q$
be positive integers. Then
the $d$-dimensional MCARMA$(p,q)$ process can be interpreted as the stationary solution to the
stochastic differential equation
\beam \label{Def: MCARMA_1}
    \mathbf{P}(D)\Ybf_t=\mathbf{Q}(D)D\Lbf_t \quad \mbox{ for } t\in\R,
\eeam
where $D$ is the differential operator, and $\Pbf$, $\Qbf$ are given as in
\eqref{Pol P} respectively \eqref{Pol Q}.
 By this representation we see the analogy
to VARMA processes: the backshift operator $B$ is replaced
by the differential operator $D$ and the iid sequence $(\xi_k)$ by the L\'{e}vy process $\Lbf$ which has independent
and stationary increments. However, this is not the formal definition
of a MCARMA process since a L\'{e}vy process is not differentiable; see \Cref{Section:MCARMA}.
The formal definition of MCARMA processes was first given in \cite{Stelzer:Marquardt2005}.
Although the history of Gaussian CARMA (the one-dimensional MCARMA) processes is very old (cf. \cite{Doob1944})
the interest in  L\'{e}vy driven CARMA processes grew quickly in the last decade;
see \cite{Brockwell2009} for an overview. The well-known
multivariate Ornstein-Uhlenbeck process is a typical example of a MCARMA process.
MCARMA processes are important for stochastic modeling in many areas of
applications as, e.g., signal processing and control (cf.
\cite{GarnierWang2008,LarssonMossbergSoederstroem2006}),
econometrics (cf.
\cite{Bergstrom1990}),  and financial
mathematics (cf. \cite{Bethetal:2014,Benthetal}).
Most of the literature restricts attention to CARMA processes which are easier to handle. Their exist only
a few references which look at MCARMA processes such as
\cite{Schlemm:Stelzer:2012}, dealing with quasi-maximum-likelihood estimation, and
\cite{Brockwell:Schlemm2013},  dealing with recovery of the driving L\'{e}vy process when the MCARMA process is sampled on a discrete-time grid.
 From \cite{Fasen11a} we already
know that the sample mean and the sample covariance of a high frequency sampled MCARMA process are consistent estimators
for the expectation respectively the covariance.

In recent years interest in the modeling of high frequency data, as they occur in finance  and turbulence,
has increased rapidly (cf. \cite{Barndorffetal:2013,Breymann:Dias:Embrecht,Todorov2009}).
  The estimation of the periodogram, normalized periodogram, smoothed periodogram
and parameter estimation in a high frequency sampled CARMA model is the topic of \cite{Fasen:Fuchs:2013a,Fasen:Fuchs:2013b}.
Moreover, \cite{Brockwell:Ferrazzano:Klueppelberg:2012b} develop a method to estimate the kernel function
of high frequency sampled MA processes, and
\cite{Ferrazzano:Fuchs} estimate the increments of the driving L\'{e}vy process in high frequency sampled MA models.

The content of this paper is the asymptotic behavior of the sample autocovariance
function of a high frequency sampled MCARMA process.
The idea is that we have data $\Ybf_{\Delta_n},\ldots,\Ybf_{n\Delta_n}$ at hand where
$\Delta_n\to0$ and $n\Delta_n\to\infty$ as $n\to\infty$.
We investigate the
asymptotic behavior of the sample autocovariance function
\begin{eqnarray} \label{ACF}
    \wh\Gamma_n(h)=\frac{1}{n}\sum_{k=1}^{n-h/\Delta_n}(\Ybf_{k\Delta_n}-\ov \Ybf_n)(\Ybf_{k\Delta_n+h}-\ov \Ybf_n)^T \quad\mbox{ for } h\in\{0,\Delta_n,\ldots,(n-1)\Delta_n\},
\end{eqnarray}
where $\ov \Ybf_n=n^{-1}\sum_{k=1}^{n}\Ybf_{k\Delta_n}$ is the sample mean, at different lags $h$. To be more precise we study
the joint asymptotic behavior of $(\wh\Gamma_n(h))_{h\in\mathcal{H}}$ for some finite set $\mathcal{H}\subseteq \bigcap_{n\geq n_0}\{0,\Delta_n,2\Delta_n,\ldots,(n-1)\Delta_n\}$, $n_0\in\N$.
We show that the  sample autocovariance function is a consistent and an asymptotically
normally distributed estimator for the autocovariance function. We present a very general
representation of the limit random matrix which helps to understand
the dependence between the components of the process quite well as, e.g., cross-covariances.
A challenge is, on the one hand, the multivariate structure of the covariances which requires
a basic knowledge of matrix calculations, and which would not be necessary if we restrict our attention only
to CARMA processes; but also for CARMA processes the results are new.
Without much effort we obtain likewise the analogous results for the sample autocovariance
function of a multivariate
MA process in discrete time  extending the work of \cite{Su:Lund}.
The structure of the limit distributions of the sample autocovariance functions of  multivariate MA and MCARMA
processes are of an analogous form. Investigating the sample autocovariance and autocorrelation function in the
one-dimensional models show that Bartlett's formula for the autocovariance and the
autocorrelation function
have the same structure
in the continuous-time and in the discrete-time model;
only sums in the discrete-time model are  integrals in the continuous-time model
and the moments of the white noise are the moments of the L\'{e}vy process.

The paper is structured on the following way.
We start with the formal definition of
a MCARMA process and a short motivation for the asymptotic
behavior of the sample autocovariance function of a high
frequency sampled MCARMA process in \Cref{Section:MCARMA}.
For the proof of the asymptotic behavior of the sample
autocovariance function we require some
preliminary limit results which are the topic of \Cref{section:limit results}. This
section is divided in two parts. The first part,
\Cref{limit continuous time}, contains limit results
for the investigation of high frequency sampled
MCARMA processes, and the second part, \Cref{Proofs Section 3},
contains the proofs.
The main section of this paper is \Cref{Section:ACF CARMA}
where we give the asymptotic behavior
of the sample autocovariance function of high
frequency sampled MCARMA processes including the asymptotic
behavior of cross-covariances between the components of a MCARMA
process and Bartlett's formula
for a CARMA process. Again a subsection
contains the proofs. All presented estimators are consistent and
asymptotically normally distributed.
Finally, in \Cref{Section:ACF MA} we introduce similar
results for  multivariate MA processes and compare both models.

\subsubsection*{Notation}
We  use the notation $\Longrightarrow$ for weak convergence and
$\stackrel{\p}{\longrightarrow}$ for convergence in probability. 
For two random vectors $\Xbf,\Ybf$ the notation $\Xbf\eqd \Ybf$ means equality in distribution.
The Euclidean norm in $\R^d$ is denoted by $ \lVert\cdot \rVert$
and the corresponding operator norm for
matrices by $ \lVert\cdot \rVert$, which is submultiplicative. Recall that two norms
on a finite-dimensional linear space are always equivalent and hence,
our results remain true if we replace the Euclidean norm by any other norm.
For $\Abf\in\R^{d\times m}$ the vec-operator $\vect(\Abf)$ is a vector in $\R^{dm}$
which is obtained by stacking the columns  of $\Abf$. The Kronecker product of two matrices
$\Abf\in\R^{d\times m},\Bbf\in\R^{l\times k}$ is denoted by
\begin{eqnarray*}
    \Abf\otimes\Bbf=\left(\begin{array}{cccc}
            A_{1,1}\Bbf & A_{1,2}\Bbf &\cdots & A_{1,m}\Bbf\\
            \vdots & \vdots & \cdots & \vdots\\
            A_{d,1}\Bbf & A_{d,2}\Bbf &\cdots & A_{d,m}\Bbf
    \end{array}\right)\in\R^{dl\times mk},
\end{eqnarray*}
where $A_{i,j}$ denotes the entry of $\Abf$ in the $i$-th row and in the $j$-th column.
 The matrix $0_{d\times m}$
is the zero matrix in $\R^{d\times m}$, $I_{d}$ is the
identity matrix in $\R^{d\times d}$ and $e_j=(0,\ldots,0,1,0,\ldots,0)\in\R^d$.
The representation $\mbox{diag}(u_1,\ldots,u_d)$
denotes a diagonal matrix in $\R^{d\times d}$ with diagonal entries $u_1,\ldots,u_d$.
For some matrix \linebreak $\Sigma\in\R^{d\times d}$ the representation $\Sigma=\Sigma^{1/2}\cdot\Sigma^{1/2 T}$
means there exists a matrix $A\in\R^{d\times d}$ such that $\Sigma=A\cdot A^T$ and $\Sigma^{1/2}:=A$.
For a vector $\xbf\in\R^d$ we write $\xbf^T$ for its
transpose and for $x\in\R$ we
write
$\lfloor x\rfloor=\sup\{k\in\Z: \,k\leq x\}$ and $\lceil x\rceil=\inf\{k\in\Z: \,x\leq k\}$. The space $(\mathbb{D}[0,T],\R^d)$
denotes  the space of all c\`{a}dl\`{a}g (continue \`{a}
droite et limit\'{e}e \`{a} gauche $=$ right continuous, with left limits)
functions  on $[0,T]$ ($T>0$) with values in $ \R^d$ equipped with the Skorokhod
$J_1$ topology.

\subsubsection*{Matrix calculation}

We would like to repeat some calculation rules for Kronecker products which are used throughout the
paper; for details we refer to \cite{Bernstein}.
Let $\xbf\in\R^n,\,\ybf\in\R^m$ be vectors and $\Abf\in\R^{n\times m},\,\Bbf\in\R^{m\times l},\,\Cbf\in\R^{l\times k}$, $\Dff\in\R^{k\times u}$
be matrices. Then
\begin{eqnarray} \label{section:matrix calculation}
\hspace*{-0.2cm}\begin{array}{r@{\hspace*{0cm}}l@{\hspace*{1cm}}r@{\hspace*{0cm}}l@{\hspace*{1cm}}c}
\xbf\ybf^T&=\xbf\otimes\ybf^T=\ybf^T\otimes\xbf, &
\vect(\xbf\ybf^T)&=\ybf\otimes \xbf,\\
\vect(\Abf\Bbf\Cbf)&=(\Cbf^T\otimes \Abf)\vect(\Bbf), &
(\Abf\otimes \Bbf)^T&=\Abf^T\otimes\Bbf^T, &
(\Abf\otimes\Cbf)(\Bbf\otimes\Dff)=(\Abf\Bbf)\otimes(\Cbf\Dff).
\end{array}
\end{eqnarray}
The matrix $P_{m,m}=\sum_{i,j=1}^me_ie_j^T\otimes e_je_i^T\in\R^{m^2\times m^2}$ is the Kronecker permutation matrix.  For $\Abf,\Bbf\in\R^{m\times m}$ and $\xbf,\ybf\in\R^m$ it has the property
\begin{eqnarray}
    P_{m,m}^2=I_{m^2}, \quad\quad
    P_{m,m}\cdot(\xbf\otimes\ybf)=\ybf\otimes\xbf,  \quad\quad \label{Permutation_matrix}
    P_{m,m}\cdot(\Abf\otimes\Bbf)=(\Bbf\otimes\Abf)P_{m,m},
\end{eqnarray}
see \cite{Bernstein}, Fact 7.4.30, where several properties of the Kronecker permutation matrix
are listed.

\section{MCARMA processes} \label{Section:MCARMA}

In this section we present some background on multivariate continuous-time ARMA (MCARMA) processes.
Since a L\'{e}vy process is not differentiable, the differential equation \eqref{Def: MCARMA_1}
cannot be used as definition of a MACARMA process.
However, it can be interpreted to be equivalent to the following definition, see~\cite{Stelzer:Marquardt2005}.

\begin{Definition} \label{Definition CARMA}
Let $(\Lbf_t)_{t\in\R}=(L_1(t),\ldots,L_m(t))_{t\in\R}$  be an $\R^m$-valued L\'{e}vy process and let
the polynomials $\mathbf{P}(z),\mathbf{Q}(z)$ be
defined as in \eqref{Pol P} and \eqref{Pol Q} with $p,q\in\N_0$, $q<p$,
and $\Qbf_0\not=0_{d\times m}$.
Moreover, define
\beao
    \Labf=-\left(\begin{array}{ccccc}
        0_{d\times d} & I_{d} & 0_{d\times d} & \cdots & 0_{d\times d}\\
        0_{d\times d} & 0_{d\times d} & I_{d} & \ddots & \vdots \\
        \vdots & & \ddots & \ddots & 0_{d\times d}\\
        0_{d\times d} & \cdots & \cdots & 0_{d\times d} & I_{d}\\
        -\Pbf_p & -\Pbf_{p-1} & \cdots & \cdots & -\Pbf_1
    \end{array}\right) \in \R^{pd\times pd},
\eeao
$\Ebf=(I_{d},0_{d\times d},\ldots,0_{d\times d}) \in \R^{d\times pd}$ and $\Bbf=(\Bbf_1^T \cdots \Bbf_p^T)^T\in \R^{pd\times m}$ with
\beao
   \Bbf_1:=\ldots:=\Bbf_{p-q-1}:=0_{d\times m}
   \quad \mbox{ and } \quad \Bbf_{p-j}:=-\sum_{i=1}^{p-j-1}\Pbf_i\Bbf_{p-j-i}+\Qbf_{q-j}  \quad \mbox{ for } j=0,\ldots,q.
\eeao
Assume $\{z\in\C:\det(\mathbf{P}(z))=0\}=\{z\in\C:\det(-\Labf- zI_{pd})=0\}\subseteq (-\infty,0)+i\R$.
Furthermore, the L\'{e}vy measure $\nu_{\Lbf}$ of $\Lbf$ satisfies
$
    \int_{\|\xbf\|>1}\log\|\xbf\|\,\nu_{\Lbf}(d\xbf)<\infty.
$
Then the $\R^{d}$-valued {\em causal MCARMA$(p,q)$ process} $(\Ybf_t)_{t\in\R}$ is  defined by the state-space
equation
\beam \label{CARMA:observation}
    \Ybf_t=\Ebf\Zbf_t \quad \mbox{ for } t\in\R,
\eeam
where
\beam \label{CARMA:state}
    \Zbf_t=\int_{-\infty}^{t}\e^{-\Labf(t-s)}\Bbf\, \dd \Lbf_s \quad \mbox{ for } t\in\R
\eeam
is the stationary unique solution to the $pd$-dimensional
stochastic differential equation \linebreak  $\dd \Zbf_t=-\Labf \Zbf_t\,\dd t+\Bbf \,\dd\Lbf_t$.
The function
$\fbf(t)=\Ebf\e^{-\Labf t}\Bbf\1_{(0,\infty)}(t)$
for $t\in\R$ is called the kernel function.
\end{Definition}

It is well known that the  stationary Ornstein-Uhlenbeck process $\Zbf$ given in \eqref{CARMA:state} observed at the
time-grid $\Delta_n\Z=\{\ldots, -2\Delta_n,-\Delta_n,0,\Delta_n,2\Delta_n,\ldots\}$ with $\Delta_n$
a positive constant has a representation as a MA process
\beao
    \Zbf_{k\Delta_n}=\sum_{j=0}^{\infty}\e^{-\Labf \Delta_nj}\xibf_{n,k-j} \quad \mbox{ for }k\in\Z,
\eeao
where $(\xi_{n,k})_{k\in\N}$ is a sequence of iid random vectors in $\R^{pd}$ with
\beam \label{xi}
  \xibf_{n,k}=\int_{(k-1)\Delta_n}^{k\Delta_n}\e^{-\Labf(k\Delta_n-s)}\, \Bbf\,\dd\Lbf_s \quad \mbox{ for } k\in\Z,n\in\N.
\eeam
To derive the asymptotic behavior of the sample autocovariance function $\wh\Gamma_n(h)$ as given in \eqref{ACF} we have to prove several intermediate steps.
First, let us define
\begin{eqnarray} \label{stern gamma}
\wh\Gamma_n^*(h)=\frac{1}{n}\sum_{k=1}^n\Ybf_{k\Delta_n}\Ybf_{k\Delta_n+h}^T.
\end{eqnarray}
 We will show
that $\sqrt{n\Delta_n}(\wh\Gamma_n^*(h)-\wh\Gamma_n(h))=o_P(1)$ so that it is sufficient to investigate
the asymptotic behavior of $\wh\Gamma_n^*(h)$. By
the Beveridge-Nelson decomposition we are able to
show that
\begin{eqnarray} \label{Short:ACF}
    \lefteqn{\sqrt{n\Delta_n}(\wh\Gamma_n^*(h)-\Gamma(h))} \nonumber\\
    &&=\Delta_n\sum_{j=0}^{\infty}\Ebf\e^{-\Labf \Delta_nj}\left(\frac{1}{\sqrt{n\Delta_n}}\sum_{k=1}^n[\xibf_{n,k}\xibf_{n,k}^T-\E(\xibf_{n,k}\xibf_{n,k}^T)]
                    \right)\e^{-\Labf^T (\Delta_nj+h)}\Ebf^T  \\
        &&\hspace*{0.3cm}+\Delta_n\sum_{j=0}^\infty\Ebf^T\e^{-\Labf\Delta_nj}\left(\frac{1}{\sqrt{n\Delta_n}}\sum_{k=1}^n\left[\sum_{r=1}^{\infty}\xibf_{n,k}\xibf_{n,k-r}^T\e^{-\Lambda^T(h+\Delta_nr)}+
            \sum_{r=1}^{\lfloor h/\Delta_n\rfloor}\xibf_{n,k}\xibf_{n,k+r}^T\e^{-\Labf^T (h-\Delta_n r)}
            \right.\right. \nonumber\\
            &&\hspace*{4cm}
            \left.\left.+\sum_{r=\lfloor h/\Delta_n\rfloor+1}^\infty\e^{-\Labf (\Delta_n r-h)}\xibf_{n,k}\xibf_{n,k+r}^T\right]\right)\e^{-\Labf^T \Delta_n j}\Ebf^T +o_P(1).
            \nonumber
\end{eqnarray}
This representation is not obvious and will first be developed on pp.~\pageref{C1}. From this we see
that we have to understand the joint limit behavior of the four terms in the brackets in \eqref{Short:ACF}, and this is what we will
do in the next section.

\section{Limit results for processes with finite fourth moments} \label{section:limit results}

\subsection{Models in continuous time} \label{limit continuous time}

The main ingredient to derive the asymptotic behavior of the sample autocovariance function
for high frequency sampled MCARMA processes is the following joint limit result of the four terms in the brackets in \eqref{Short:ACF}.

\begin{Proposition} \label{Theorem: limit MCARMA}
 Let $(\Ybf_t)_{t\geq 0}$ be a MCARMA process as defined in \eqref{CARMA:observation} with $\E\|\Lbf_1\|^4<\infty$,
 $\E(\Lbf_1)=0_m$
 and $\E(\Lbf_1\Lbf_1^T)=\Sigma=\Sigma^{1/2}\cdot\Sigma^{1/2\,T}\in\R^{m\times m}$.  The sequence
 $(\xi_{n,k})$ is defined as in \eqref{xi} and
\begin{eqnarray*}
    \Upsilon:=\int_{\R^m}\xbf\xbf^T\otimes \xbf\xbf^T\,\nu(d\xbf).
\end{eqnarray*}
Suppose $(\Delta_n)_{n\in\N}$ is a sequence of positive constants
with $\Delta_n\downarrow 0$ and $\lim_{n\to\infty}n\Delta_n=\infty$.  We assume there exists a sequence of positive constants
$l_n\to\infty$ with  $n/l_n\to\infty$ and $l_n\Delta_n\to\infty$.
Let $\mathcal{H}\subseteq\left[0,\infty\right)$ be a finite set. Then as $n\to\infty$,
\begin{eqnarray*}
    \lefteqn{\hspace*{-1cm}\left(\frac{1}{\sqrt{n\Delta_n}}\sum_{r=1}^{\infty}\sum_{k=1}^n[\xibf_{n,k}\xibf_{n,k-r}^T]\e^{-\Lambda^T(h+\Delta_nr)},
            \frac{1}{\sqrt{n\Delta_n}}\sum_{r=1}^{\lfloor h/\Delta_n\rfloor}\sum_{k=1}^n\xibf_{n,k}\xibf_{n,k+r}^T\e^{-\Labf^T (h-\Delta_n r)},\right.}\\
            &&\left.\hspace*{-0.5cm}\frac{1}{\sqrt{n\Delta_n}}\sum_{r=\lfloor h/\Delta_n\rfloor+1}^{\infty}\e^{-\Labf (\Delta_n r-h)}\sum_{k=1}^n[\xibf_{n,k}\xibf_{n,k+r}^T],
            \frac{1}{\sqrt{n \Delta_n}}\sum_{k=1}^n[\xibf_{n,k}\xibf_{n,k}^T-\E(\xibf_{n,1}\xibf_{n,1}^T)]
            \right)_{h\in \mathcal{H}}\\
    &&\hspace*{-1cm}\stackrel{\mathcal{D}}{\to}\left(\int_0^\infty\Bbf\Sigma_{\Lbf}^{1/2}\,{\rm}d\Wbf_s\Sigma^{1/2\,T}_{\Lbf}\Bbf^T\e^{-\Lambda^T(h+s)},
    \int_0^h\Bbf\Sigma_{\Lbf}^{1/2}\,{\rm}d\Wbf_s^T\Sigma^{1/2\,T}_{\Lbf}\Bbf^T\e^{-\Lambda^T(h-s)},\right.\\
    &&\hspace*{-0.5cm}\left.\int_h^\infty\e^{-\Lambda (s-h)}\Bbf\Sigma^{1/2}_{\Lbf}\,{\rm}d\Wbf_s^T\Sigma^{1/2\,T}_{\Lbf}\Bbf^T,\Bbf\Wbf^*(\Upsilon)\Bbf^T\right)_{h\in\mathcal{H}},
\end{eqnarray*}
where $\Wbf^*(\Upsilon)$
 is an $\R^{pd\times pd}$-valued normal random matrix with $\vect(\Wbf^*(\Upsilon))\sim\mathcal{N}(0_{(pd)^2},\Upsilon)$
 independent from the $\R^{pd\times pd}$-valued standard Brownian motion $(\Wbf_t)_{t\geq 0}$.
\end{Proposition}

\begin{Remark}  \label{Remark_1} $\mbox{}$
We investigate in detail the convergence of the last term.
\begin{itemize}
\item[(a)] Let $\Lbf$ be a Brownian motion. Then $\nu=0$ and hence, $\Upsilon=0_{m^2\times m^2}$. Thus, a conclusion of \cref{Theorem: limit MCARMA}
is that as $n\to\infty$,
\begin{eqnarray*}
    \frac{1}{\sqrt{n \Delta_n}}\sum_{k=1}^n[\xibf_{n,k}\xibf_{n,k}^T-\E(\xibf_{n,1}\xibf_{n,1}^T)]\stackrel{\mathcal{D}}{\to}0_{(pd)^2\times (pd)^2}.
\end{eqnarray*}
\item[(b)]
When $\Lbf$ has independent components then $\Wbf^*(\Upsilon)$ reduces
to a much simpler random matrix. Define $\theta_i:=\int_{\R} x^4\,\nu_i(d x)=\E(L_i(1)^4)-3\E(L_i(1)^2)$, $i=1,\ldots,m$, where
$L_i$ is the $i$-th component of $\Lbf$ with L\'{e}vy measure $\nu_i$. Then
    $\Upsilon= \sum_{i=1}^m\left[(e_ie_i^T\otimes e_ie_i^T)\cdot\theta_i\right],$
and thus,  $$\Wbf^*(\Upsilon)\eqd \sum_{i=1}^m e_i\sqrt{\theta_i}N_i e_i^T
= \diag(\sqrt{\theta_1}N_1,\ldots,\sqrt{\theta_m}N_m), $$ where $N_1,\ldots,N_m$
are iid $\mathcal{N}(0,1)$-distributed.
In particular, we obtain in the one-dimensional case as $n\to\infty$,
\begin{eqnarray} \label{NA1}
    \frac{1}{\sqrt{n \Delta_n}}\sum_{k=1}^n[\xi_{n,k}^2-\E(\xi_{n,1}^2)]\weak\mathcal{N}(0,\E(L_1^4)-3\E(L_1^2)).
\end{eqnarray}
However, if $(\xi_k)_{k\in\N}$ is an iid sequence with $\E(\xi_k)=0$ and $\E|\xi_k|^4<\infty$ then
obviously by the classical central limit theorem of Lindeberg-L\'{e}vy as $n\to\infty$,
\begin{eqnarray} \label{NA2}
     \frac{1}{\sqrt{n }}\sum_{k=1}^n[\xi_{k}^2-\E(\xi_{1}^2)]\weak\mathcal{N}(0,\E(\xi_1^4)-\E(\xi_1^2)^2).
\end{eqnarray}
The limit in \eqref{NA1} has a factor 3 which does not appear in \eqref{NA2}.
\hfill$\Box$
\end{itemize}
\end{Remark}

The proof of \cref{Theorem: limit MCARMA} is based on some limit results which are interesting on their own. The main task is
to derive \cref{CLT_2}.

\begin{Proposition} \label{CLT_2}
Let $(\xi_{n,k})_{k\in\bbn}$ be a sequence of iid random vectors in $\mathbb{R}^{pd}$
with $\E\|\xi_{n,k}\|^4<\infty$, \linebreak $\E(\xi_{n,1})=0_{pd}$ and $\E(\xi_{n,1}\xi_{n,1}^T)=\Sigma_n=\Sigma_n^{1/2}\cdot \Sigma_n^{1/2 T}\in \mathbb{R}^{pd\times pd}$
 for any $n\in\bbn$.
Suppose $(\Delta_n)_{n\in\N}$ is a sequence of positive constants
with $\Delta_n\downarrow 0$ and $\lim_{n\to\infty}n\Delta_n=\infty$.  We assume that there exists a sequence of positive constants
$l_n\to\infty$ with  $n/l_n\to\infty$ and $l_n\Delta_n\to\infty$. Moreover,
$\Delta_n^{-1}\Sigma_n\to\Sigma=\Sigma^{1/2}\cdot \Sigma^{1/2 T}\in \mathbb{R}^{pd\times pd}$ as $n\to\infty$,
$\E\|\xi_{n,1}\|^4\leq \mbox{const.}\cdot\Delta_n$\, for any $n\in\N$ and
\begin{eqnarray} \label{A.59}
    \lim_{n\to\infty}\Delta_n^{-1}\E(\xi_{n,1}\xi_{n,1}^T\otimes \xi_{n,1}\xi_{n,1}^T)=\Upsilon\in\R^{(pd)^2\times(pd)^2}.
\end{eqnarray}
Define for $t\geq 0$, $n\in\N$,
\begin{eqnarray*}
    S_n^{(1)}(t)=\sum_{j=1}^{\lfloor t/\Delta_n \rfloor}\sum_{k=1}^n\xi_{n,k}\xi_{n,k-j}^T, \quad
    S_n^{(2)}(t)=\sum_{j=1}^{\lfloor t/\Delta_n \rfloor}\sum_{k=1}^n\xi_{n,k}\xi_{n,k+j}^T
    \quad \mbox{ and } \quad
    S_n^{(3)}=\sum_{k=1}^n[ \xibf_{n,k}\xibf_{n,k}^T-\Sigma_{n}].
\end{eqnarray*}
Then as $n\to\infty$,
\begin{eqnarray*}
     \frac{1}{\sqrt{n\Delta_n}}(S_n(t))_{t\in[0,T]}&:=&{\left(\frac{1}{\sqrt{n\Delta_n}}S_n^{(1)}(t),\frac{1}{\sqrt{n\Delta_n}}S_n^{(2)}(t),\frac{1}{\sqrt{n\Delta_n}}S_n^{(3)}\right)_{t\in[0,T]}}\\
     &\stackrel{\mathcal{D}}{\to}&(\Sigma^{1/2}\mathbf{W}_t\Sigma^{1/2\,T},\Sigma^{1/2}\mathbf{W}_t^T\Sigma^{1/2\,T},\Wbf^*(\Upsilon))_{t\in[0,T]}
\end{eqnarray*}
 in $D(\left[0,T\right],\mathbb{R}^{pd\times pd}\times\mathbb{R}^{pd\times pd}\times\mathbb{R}^{pd\times pd})$ where $\Wbf^*(\Upsilon)$
 is an $\R^{pd\times pd}$-valued normal random matrix with $\vect(\Wbf^*(\Upsilon))\sim\mathcal{N}(0_{(pd)^2},\Upsilon)$
 independent from the $\R^{pd\times pd}$-valued standard Brownian motion $(\Wbf_t)_{t\geq 0}$.
\end{Proposition}

A conclusion of \cref{CLT_2} and a continuous mapping theorem is \cref{Lemma1.3}.
\cref{Theorem: limit MCARMA} can be seen as special case of  \cref{Lemma1.3}, we have only
to check that the assumptions are satisfied.

\begin{Proposition} \label{Lemma1.3}
 Let the assumptions of \cref{CLT_2} hold. Suppose $\gbf_l:\left[0,\infty\right)\to\R^{pd\times pd}$  ($l=1,\ldots,M$) are maps
 with finite variation and $\int_0^\infty\|\gbf_l(s)\|^2\,ds<\infty$. Then as $n\to\infty$,
\begin{eqnarray*}
    \lefteqn{\frac{1}{\sqrt{n\Delta_n}}\left(\int_0^{\infty}\, \gbf_l(s)S_n^{(1)}(\mbox{d}s),\int_0^{\infty}\, S_n^{(1)}(\mbox{d}s)\gbf_l(s),\int_0^{\infty} \gbf_l(s)S_n^{(2)}(\mbox{d}s),\int_0^{\infty} S_n^{(2)}(\mbox{d}s)\gbf_l(s)
    ,S_n^{(3)}\right)_{l=1,\ldots,M}}\\
    &&\stackrel{\mathcal{D}}{\to}\left(\int_0^\infty\gbf_l(s)\Sigma^{1/2}\,{\rm}d\Wbf_s\Sigma^{1/2\,T},\int_0^\infty\Sigma^{1/2}\,{\rm}d\Wbf_s\Sigma^{1/2\,T}\gbf_l(s),
    \int_0^\infty\gbf_l(s)\Sigma^{1/2}\,{\rm}d\Wbf_s^T\Sigma^{1/2\,T},\right.\\
    &&\quad\quad\left.\int_0^\infty\Sigma^{1/2}\,{\rm}d\Wbf_s^T\Sigma^{1/2\,T}\gbf_l(s)^T,
    \Wbf^*(\Upsilon)\right)_{l=1,\ldots,M}.
\end{eqnarray*}

\end{Proposition}

Now for the proof of \cref{Theorem: limit MCARMA} we have mainly to check that the assumptions of \cref{Lemma1.3}
are satisfied, in particularly that \eqref{A.59} holds.
\cite{Asmussen:Rosinski}, Lemma~3.1, already derived the limit behavior
$\lim_{n\to\infty}\Delta_n^{-1}\E( L_{\Delta_n}^4)$ for an one-dimensional L\'{e}vy process. We have to
extend this result to a multivariate L\'{e}vy process and use it to show \eqref{A.59}.

\begin{Lemma}  \label{Appendix:1} $\mbox{}$
\begin{itemize}
\item[(a)] Let $\Lbf=(\Lbf_t)_{t\geq 0}$ be a L\'{e}vy process with $\E\|\Lbf_1\|^4<\infty$, $\E(\Lbf_1)=0_m$, and $(\Delta_n)_{n\in\N}$
be a sequence of positive constants with $\lim_{n\to\infty}\Delta_n=0$. Then
    $${\displaystyle \lim_{n\to\infty}\Delta_n^{-1}\E(\Lbf_{\Delta_n}\Lbf_{\Delta_n}^T\otimes \Lbf_{\Delta_n}\Lbf_{\Delta_n}^T)
        =\int_{\R^m} \xbf\xbf^T\otimes \xbf\xbf^T\,\nu(d\xbf)}.$$
\item[(b)]  Let the assumptions of \cref{Theorem: limit MCARMA} hold. Then
    $${\displaystyle \lim_{n\to\infty}\Delta_n^{-1}\E(\xi_{n,1}\xi_{n,1}^T\otimes \xi_{n,1}\xi_{n,1}^T)
        =\Bbf\otimes\Bbf\left(\int_{\R^m} \xbf\xbf^T\otimes \xbf\xbf^T\,\nu(d\xbf)\right)\Bbf^T\otimes\Bbf^T.}$$
\end{itemize}
\end{Lemma}

\subsection{Proofs} \label{Proofs Section 3}

\subsubsection{Auxiliary results for the proof of \cref{CLT_2}}

For the proof of \cref{CLT_2} we derive some auxiliary results. First, we want to characterize
the limit process $(\Sigma^{1/2}\Wbf_t\Sigma^{1/2 T},\Sigma^{1/2}\Wbf_t^T\Sigma^{1/2 T})_{t\geq 0}$.

\begin{Lemma} \label{limits_structure}
Let $P_{pd,pd}=\sum_{i,j=1}^{pd}e_ie_j^T\otimes e_je_i^T\in\R^{(pd)^2\times (pd)^2}$ be the Kronecker
permutation matrix and
    \begin{eqnarray} \label{SigmaStar}
        \Sigma^*:=
        \left(\begin{array}{cc}
            \Sigma^{1/2}\otimes \Sigma^{1/2} & 0_{(pd)^2\times (pd)^2}\\
            0_{(pd)^2\times (pd)^2}   & \Sigma^{1/2}\otimes \Sigma^{1/2}
        \end{array}\right)\cdot
        \left(\begin{array}{cc}
            I_{(pd)^2} &  P_{pd,pd}\\
            P_{pd,pd} & I_{(pd)^2}
        \end{array}\right).
        \left(\begin{array}{cc}
            \Sigma^{1/2}\otimes \Sigma^{1/2} & 0_{(pd)^2\times (pd)^2}\\
            0_{(pd)^2\times (pd)^2}   & \Sigma^{1/2}\otimes \Sigma^{1/2}
        \end{array}\right)^T.
    \end{eqnarray}
    Then
    $(\vect(\Sigma^{1/2}\Wbf_t\Sigma^{1/2 T},\Sigma^{1/2}\Wbf_t^T\Sigma^{1/2 T}))_{t\geq 0}$
    is a Brownian motion with covariance matrix $\Sigma^*$.
\end{Lemma}
\paragraph{Proof.}
    The reason is that since $\vect(\Sigma^{1/2}\Wbf_1\Sigma^{1/2 T})=\Sigma^{1/2}\otimes \Sigma^{1/2}\vect(\Wbf_1)$
    we have
    \begin{eqnarray*}
    \vect(\Sigma^{1/2}\Wbf_1\Sigma^{1/2 T},\Sigma^{1/2}\Wbf_1^T\Sigma^{1/2 T})=\left(\begin{array}{cc}
            \Sigma^{1/2}\otimes \Sigma^{1/2} & 0_{(pd)^2\times (pd)^2}\\
            0_{(pd)^2\times (pd)^2}   & \Sigma^{1/2}\otimes \Sigma^{1/2}
        \end{array}\right)\vect(\Wbf_1,\Wbf_1^T).
    \end{eqnarray*}
    Some straightforward calculations give
     \begin{eqnarray} \label{moment Kronecker}
         \E(\vect(\Wbf_1)\vect(\Wbf_1^T)^T)=\E(\Wbf_1\otimes\Wbf_1^T)=\sum_{i,j=1}^{pd}e_ie_j^T\otimes e_je_i^T=P_{pd,pd},
    \end{eqnarray}
    and thus,
    \begin{eqnarray*}
        \E(\vect(\Sigma^{1/2}\Wbf_1\Sigma^{1/2 T},\Sigma^{1/2}\Wbf_1^T\Sigma^{1/2 T})\cdot \vect(\Sigma^{1/2}\Wbf_1\Sigma^{1/2 T},\Sigma^{1/2}\Wbf_1^T\Sigma^{1/2 T})^T)
        =\Sigma^*.
    \end{eqnarray*}
    The stationary and independent increment property of the Brownian motion $(\Wbf_t)_{t\geq 0}$
    transfer to \linebreak $(\vect(\Sigma^{1/2}\Wbf_t\Sigma^{1/2 T},\Sigma^{1/2}\Wbf_t^T\Sigma^{1/2 T}))_{t\geq 0}$
    such that the conclusion follows.
\hfill$\Box$\\

Next we prove the convergence of $S_n^{(3)}$ alone which is more or less straightforward.

\begin{Lemma} \label{A.599}

   ${\displaystyle     \frac{1}{\sqrt{n\Delta_n}}\sum_{k=1}^n[\xi_{n,k}\otimes\xi_{n,k}-\vect(\Sigma_n)]\weak \mathcal{N}(0_{(pd)^2},\Upsilon)}$
   as $n\to\infty$.

\end{Lemma}
\paragraph{Proof.}
    By assumption  $(\xi_{n,k}\otimes  \xi_{n,k})_{k\in\N}$ is a sequence of iid random vectors with
    $\E(\xi_{n,k}\otimes  \xi_{n,k})=\vect({\Sigma}_n)$ and
    \begin{eqnarray*}
        \lim_{n\to\infty}\Delta_n^{-1}[\E((\xi_{n,k}\otimes  \xi_{n,k})\cdot(\xi_{n,k}\otimes  \xi_{n,k})^T)-\vect({\Sigma}_n)\vect({\Sigma}_n)^T]=\Upsilon,
    \end{eqnarray*}
    where we used that $\lim_{n\to\infty}\Delta_n^{-1}\vect({\Sigma}_n)\vect({\Sigma}_n)^T=0_{(pd)^2\times(pd)^2}$.
    It remains to show the Lindeberg-condition so that we can apply the central limit theorem
    of Lindeberg-Feller (see \cite{jacodB}, Theorem VII.5.2). Let $\epsilon>0$. As in \cite[Proposition A.1(d)]{Fasen11a}
    using \cite[Lemma~3.1]{Asmussen:Rosinski} it is possible to prove that
    \begin{eqnarray*}
        \lim_{n\to\infty} \Delta_n^{-1}\E(\|\xi_{n,1}\|^4\1_{\{\|\xi_{n,1}\|>\epsilon\sqrt{n\Delta_n}\}})=0.
    \end{eqnarray*}
    Thus, the central limit theorem of Lindeberg-Feller 
    gives the desired weak convergence as $n\to\infty$.
\hfill$\Box$\\

Now we are able to prove the convergence of the two-dimensional distribution in \cref{CLT_2} before we prove the convergence of the stochastic process.

\begin{Lemma} \label{Convergence:fidis}
Let $0\leq s<t<\infty$. Then as $n\to\infty$,
   \begin{eqnarray*}
        \lefteqn{\frac{1}{\sqrt{n\Delta_n}}(S_n^{(1)}(s),S_n^{(1)}(t)-S_n^{(1)}(s),S_n^{(2)}(s),S_n^{(2)}(t)-S_n^{(2)}(s),S_n^{(3)})}\\
        &&\weak
        (\Sigma^{1/2}\mathbf{W}_s\Sigma^{1/2\,T},\Sigma^{1/2}(\mathbf{W}_t-\mathbf{W}_s)\Sigma^{1/2\,T},\Sigma^{1/2}\mathbf{W}_s^T\Sigma^{1/2\,T},\Sigma^{1/2}
        (\mathbf{W}_t^T-\mathbf{W}_s^T)\Sigma^{1/2\,T},\Wbf^*(\Upsilon)) .
    \end{eqnarray*}
\end{Lemma}
\paragraph{Proof.}
The proof uses Cram\'{e}r-Wold theorem. Thus, let
$c_1,c_2\in\R^{2(pd)^2}$, $c_3\in \R^{(pd)^2}$ and define
\begin{eqnarray} \label{sn*}
    S_n^*&:=&c_1^T\vect(S_n^{(1)}(s),S_n^{(2)}(s))+c_2^T\vect(S_n^{(1)}(t)-S_n^{(1)}(s),S_n^{(2)}(t)-S_n^{(2)}(s))+c_3^T\vect(S_n^{(3)}) \nonumber\\
        &=&\sum_{k=1}^n\left[\sum_{j=1}^{\lfloor s/\Delta_n \rfloor}c_1^T\vect(\xi_{n,k}\vect(\xi_{n,k-j},\xi_{n,k+j})^T)
            +\sum_{j=\lfloor s/\Delta_n \rfloor +1}^{\lfloor t/\Delta_n \rfloor}c_2^T\vect(\xi_{n,k}\vect(\xi_{n,k-j},\xi_{n,k+j})^T)\right. \nonumber\\
            &&\left.+c_3^T\vect(\xi_{n,k}\xi_{n,k}^T-\Sigma_n)\right] \nonumber\\
        &=:&\sum_{k=1}^nZ_{n,k},
\end{eqnarray}
with
\begin{eqnarray*}
Z_{n,k}&:=&\sum_{j=1}^{\lfloor s/\Delta_n \rfloor}c_1^T\vect(\xi_{n,k}\vect(\xi_{n,k-j},\xi_{n,k+j})^T)
            +\sum_{j=\lfloor s/\Delta_n \rfloor +1}^{\lfloor t/\Delta_n \rfloor}c_2^T\vect(\xi_{n,k}\vect(\xi_{n,k-j},\xi_{n,k+j})^T)\\
            &&\quad+c_3^T\vect(\xi_{n,k} \xi_{n,k}^T-\Sigma_n).
\end{eqnarray*}
Moreover, define $Z_{n,k,j}:=\vect(\xi_{n,k}\vect(\xi_{n,k-j},\xi_{n,k+j})^T)$
    such that
    \begin{eqnarray} \label{NN1}
        Z_{n,k}=\sum_{j=1}^{\lfloor s/\Delta_n \rfloor}c_1^TZ_{n,k,j}
            +\sum_{j=\lfloor s/\Delta_n \rfloor +1}^{\lfloor t/\Delta_n \rfloor}c_2^TZ_{n,k,j}+c_3^T\vect(\xi_{n,k} \xi_{n,k}^T-\Sigma_n).
    \end{eqnarray}
    We will prove that $\frac{1}{\sqrt{n\Delta_n}}S_n^*$ converges weakly to a normal distribution with mean $0$ and variance
    \begin{eqnarray} \label{eq:variance}
        sc_1^T\Sigma^*c_1+(t-s)c_2^T\Sigma^*c_2+c_3^T\Upsilon c_3=:\Sigma(c_1,c_2,c_3),
    \end{eqnarray}
    where $\Sigma^*$ is as in \eqref{SigmaStar}.
    We take the sequence $(l_n)$ where $l_n\to\infty$,  $n/l_n\to\infty$ and $l_n\Delta_n\to\infty$
    and assume for the ease of notation that $l_n$ and $t\Delta_n^{-1}$ are integers. Write
    \begin{eqnarray} \label{sn*2}
        S_n^*&=&\left[\sum_{i=1}^{\lfloor n/l_n\rfloor}\sum_{k=(i-1)l_n+t\Delta_n^{-1}+1}^{il_n-t\Delta_n^{-1}-1}Z_{n,k}\right]
                +\left[\sum_{k=1}^{t\Delta_n^{-1}}Z_{n,k}+
                    \sum_{i=1}^{\lfloor n/l_n\rfloor-1}\sum_{k=il_n-t\Delta_n^{-1}}^{il_n+t\Delta_n^{-1}}Z_{n,k}
                +\sum_{k=\lfloor n/l_n\rfloor l_n-t\Delta_n^{-1}}^{n}Z_{n,k}\right] \nonumber\\
            &=:&S_{1,n}^*+S_{2,n}^*.
    \end{eqnarray}
    The proof is divided in two parts. On the one hand, we have to show that \linebreak $\frac{1}{\sqrt{n\Delta_n}}S_{1,n}^*\stackrel{\mathcal{D}}{\to}\mathcal{N}(0,\Sigma(c_1,c_2,c_3))$
    and on the other hand, that $\frac{1}{\sqrt{n\Delta_n}}S_{2,n}^*\stackrel{\p}{\to}0$.
    We derive the weak convergence of the first term with the central limit theorem of Lindeberg-Feller.
    Therefore, we require some auxiliary results: the asymptotic behavior of the covariance matrix
    of $S_{1,n}^*$ (\cref{Cov S}) and the Lindeberg condition (\cref{Lemma:Lindeberg}).
    \begin{Lemma} \label{Cov S}
    Let $\Sigma(c_1,c_2,c_3)$ be given as in \eqref{eq:variance} and $S_{1,n}^*$ as in \eqref{sn*2}. Then
    \begin{eqnarray*}
        \lim_{n\to\infty}(n\Delta_n)^{-1}\E(S_{1,n}^*S_{1,n}^{*\, T})=\Sigma(c_1,c_2,c_3).
    \end{eqnarray*}
    \end{Lemma}
    \paragraph{Proof.}
    We start by calculating the
    asymptotic covariance matrix of
    \begin{eqnarray*}
        Z_{n,k,j}=\left(\xi_{n,k-j}\atop \xi_{n,k+j}\right)\otimes \xi_{n,k} \quad \mbox{ where } \quad
        Z_{n,k,j} Z_{n,m,l}^{T}=\left(\begin{array}{cc}\xi_{n,k-j}\xi_{n,m-l}^T\quad & \xi_{n,k-j}\xi_{n,m+l}^T\\
        \xi_{n,k+j}\xi_{n,m-l}^T\quad & \xi_{n,k+j}\xi_{n,m+l}^T
        \end{array}
        \right)\otimes \xi_{n,k} \xi_{n,m}^T.
    \end{eqnarray*}
    Having in mind that $(\xi_{n,k})_{k\in\N}$ is an iid sequence with $\E(\xi_{n,k})=0_{pd}$, we get on the one hand,
    \begin{eqnarray*}
        \Delta_n^{-2}\E(Z_{n,k,j} Z_{n,k,j}^{T})
        \stackrel{n\to\infty}{\to}
        \left(\begin{array}{cc}(\Sigma^{1/2}\otimes \Sigma^{1/2})(\Sigma^{1/2}\otimes \Sigma^{1/2})^T \quad & 0_{(pd)^2\times (pd)^2}\\
        0_{(pd)^2\times (pd)^2}\quad & (\Sigma^{1/2}\otimes \Sigma^{1/2})(\Sigma^{1/2}\otimes \Sigma^{1/2})^T \end{array}
        \right)=:\Sigma^*_1,
    \end{eqnarray*}
    and on the other hand,
    \begin{eqnarray*}
        \Delta_n^{-2}\E(Z_{n,k,j} Z_{n,k+j,j}^{T})&=&\Delta_n^{-2}\left(\begin{array}{cc}
                0_{(pd)^2\times (pd)^2}\quad & 0_{(pd)^2\times (pd)^2} \\
                \E((\xi_{n,k+j}\xi_{n,k}^T)\otimes (\xi_{n,k+j}\xi_{n,k}^T)^T)\quad & 0_{(pd)^2\times (pd)^2}
        \end{array}
        \right)\\
        &\stackrel{n\to\infty}{\to}&\left(\begin{array}{cc}
                0_{(pd)^2\times (pd)^2} \quad & 0_{(pd)^2\times (pd)^2} \\
                (\Sigma^{1/2}\otimes \Sigma^{1/2})P_{pd,pd}(\Sigma^{1/2}\otimes \Sigma^{1/2})^T\quad & 0_{(pd)^2\times (pd)^2}
        \end{array}
        \right)=:\Sigma^*_2,
    \end{eqnarray*}
    compare \eqref{moment Kronecker}.
    Moreover, $\E(Z_{n,k,j} Z_{n,m,l}^{T})=0_{2(pd)^2\times 2(pd)^2}$ if  $(m,l)\notin\{(k,j),(k+j,j),(k-j,j)\}$.
Finally, by \eqref{NN1} and \eqref{sn*2} we obtain with
    \begin{eqnarray*}
        \lefteqn{(n\Delta_n)^{-1}\E(S_{1,n}^*S_{1,n}^{*T})} \nonumber \\
         &=&(n\Delta_n)^{-1}\lfloor n/l_n \rfloor \sum_{k=t\Delta_n^{-1}+1}^{l_n-t\Delta_n^{-1}}\left[\sum_{j=1}^{\lfloor s/\Delta_n\rfloor}c_1^T
         \E(Z_{n,k,j}Z_{n,k,j}^T)c_1+\sum_{j=\lfloor s/\Delta_n\rfloor+1}^{\lfloor t/\Delta_n\rfloor}c_2^T\E(Z_{n,k,j}Z_{n,k,j}^T)c_2  \right.\nonumber\\
         &&\quad\quad\quad\left.+\sum_{j=1}^{\lfloor s/\Delta_n\rfloor}c_1^T\left[\E(Z_{n,k,j}Z_{n,k+j,j}^T)\1_{\{k+j\leq l_n-t\Delta_n^{-1}\}}+\E(Z_{n,k,j}Z_{n,k-j,j}^T)\1_{\{k-j\geq t\Delta_n^{-1}+1\}}\right]c_1\right.\nonumber\\
         &&\left.\quad\quad\quad+\sum_{j=\lfloor s/\Delta_n\rfloor+1}^{\lfloor t/\Delta_n\rfloor}c_2^T\left[\E(Z_{n,k,j}Z_{n,k+j,j}^T)\1_{\{k+j\leq l_n-t\Delta_n^{-1}\}}+\E(Z_{n,k,j}Z_{n,k-j,j}^T)\1_{\{k-j\geq t\Delta_n^{-1}+1\}}\right]c_2\right.\nonumber\\
         &&\left.\quad\quad\quad +c_3^T\E(\vect(\xi_{n,k} \xi_{n,k}^T-\Sigma_n)\vect(\xi_{n,k} \xi_{n,k}^T-\Sigma_n)^T)c_3
         \right] \nonumber\\
            &\stackrel{n\to\infty}{\to}&sc_1^T[\Sigma_1^*+\Sigma_2^*+\Sigma_2^{*T}]c_1+(t-s)c_2^T[\Sigma_1^*+\Sigma_2^*+\Sigma_2^{*T}]c_2+c_3^T\Upsilon c_3
            =\Sigma(c_1,c_2,c_3) \label{A.111}
    \end{eqnarray*}
    the desired result.
     \hfill$\Box$

    \begin{Lemma} \label{Lemma:Lindeberg}
    The Lindeberg condition
    \begin{eqnarray*}
        \lim_{n\to\infty}\frac{1}{n\Delta_n}\sum_{i=1}^{\lfloor n/l_n\rfloor}\E\left(\sum_{k=(i-1)l_n+t\Delta_n^{-1}+1}^{il_n-t\Delta_n^{-1}-1}Z_{n,k} \1_{\left\{\left|\sum_{k=(i-1)l_n+t\Delta_n^{-1}+1}^{il_n-t\Delta_n^{-1}-1}Z_{n,k}\right|> \epsilon\sqrt{n\Delta_n}\right\}}\right)^2 =0, \quad \epsilon>0,
    \end{eqnarray*}
    is satisfied.
    \end{Lemma}
     \paragraph{Proof.}
     In connection to \eqref{NN1} let us define for $n,i\in\N$,
     \begin{eqnarray*}
        \wt Z_{n,i}=\sum_{k=(i-1)l_n+t\Delta_n^{-1}+1}^{il_n-t\Delta_n^{-1}-1}c_3^T\vect(\xi_{n,k} \xi_{n,k}^T-\Sigma_n) \quad \mbox{ and } \quad
        Z_{n,i}^*=\left(\sum_{k=(i-1)l_n+t\Delta_n^{-1}+1}^{il_n-t\Delta_n^{-1}-1}Z_{n,k}\right)-\wt Z_{n,i}.
     \end{eqnarray*}
     Then
     \begin{eqnarray}\label{A.113}
        \lefteqn{\hspace*{-1cm}\E\left(\sum_{k=(i-1)l_n+t\Delta_n^{-1}+1}^{il_n-t\Delta_n^{-1}-1}Z_{n,k} \1_{\left\{\left|\sum_{k=(i-1)l_n+t\Delta_n^{-1}+1}^{il_n-t\Delta_n^{-1}-1}Z_{n,k}\right|> \epsilon\sqrt{n\Delta_n}\right\}}\right)^2} \nonumber\\
        &&\hspace*{-1cm}\leq 2[\E(|Z_{n,1}^*|^2\1_{\{|Z_{n,1}^*|>\epsilon/2 \sqrt{n\Delta_n}\}})+
            \E(|\wt Z_{n,1}|^2\1_{\{|\wt Z_{n,1}|>\epsilon/2 \sqrt{n\Delta_n}\}})].
    \end{eqnarray}
    By the central limit result in~\cref{A.599}, $l_n/n\to 0$ as $n\to\infty$ and
    \cite[Example 28.4]{billingsley:1986} the Lindeberg-condition
    \begin{eqnarray} \label{A.57}
        \lim_{n\to\infty}\frac{1}{l_n\Delta_n}\E(|\wt Z_{n,1}|^2\1_{\{|\wt Z_{n,1}|>\epsilon/2 \sqrt{n\Delta_n}\}})=0
    \end{eqnarray}
    holds. For the second term in \eqref{A.113} we use the Ljapunov condition. Therefore, note that \linebreak
     $   \frac{1}{(l_n\Delta_n)^2}\E\|Z_{n,1}^*\|^4\leq \mbox{const.}$
    To see this,
    define the vectors $c_i^{(1)}=(c_{i,1},\ldots,c_{i,(pd)^2})\in\R^{(pd)^2}$ which contains
    the first $(pd)^2$-components of the vector $c_i$ respectively,
    $c_i^{(2)}=(c_{i,(pd)^2+1},\ldots,c_{i,2(pd)^2})\in\R^{(pd)^2}$ which contains the
    last $(pd)^2$-components of the vector $c_i$  ($i=1,2$). Moreover, for $n,k\in\N$,
    \begin{eqnarray} \label{def Z}
        Z_{n,k}^{(1)}&:=&\sum_{j=1}^{\lfloor s/\Delta_n \rfloor}c_1^{(1)\,T}\cdot (\xi_{n,k}\otimes \xi_{n,k-j}^T), \quad\quad\quad\quad\quad
                Z_{n,k}^{(2)}:=\sum_{j=1}^{\lfloor s/\Delta_n \rfloor}c_1^{(2)\,T}\cdot( \xi_{n,k}\otimes \xi_{n,k+j}^T), \\
        Z_{n,k}^{(3)}&:=&\sum_{j=\lfloor s/\Delta_n \rfloor +1}^{\lfloor t/\Delta_n \rfloor}c_2^{(1)\,T} \cdot(\xi_{n,k}\otimes \xi_{n,k-j}^T), \quad\quad\quad
        Z_{n,k}^{(4)}:=\sum_{j=\lfloor s/\Delta_n \rfloor +1}^{\lfloor t/\Delta_n \rfloor}c_2^{(2)\,T} \cdot(\xi_{n,k}\otimes \xi_{n,k+j}^T). \nonumber
    \end{eqnarray}
    Since by assumption $(\xi_{n,k})_{k\in\N}$ is an iid sequence with $\E(\xi_{n,k})=0_{pd}$, $\E\|\xi_{n,k}\|^2<\mbox{const. }\cdot\Delta_n$ and \linebreak
    $\E\|\xi_{n,k}\|^4<\mbox{const. }\cdot\Delta_n$, the sequence of random variables
    $(Z_{n,k}^{(i)})_{k\in\N}$  is an uncorrelated sequence with $\E(Z_{n,k}^{(i)})=0$, $\E((Z_{n,k}^{(i)})^2)\leq \mbox{const. }\cdot \Delta_n$
    and $\E((Z_{n,k}^{(i)})^4)\leq \mbox{const. }\cdot \Delta_n$ ($i=1,\ldots,4$). Thus,
    \begin{eqnarray*}
        \E((Z_{n,1}^*)^4)\leq\mbox{const.}\sum_{i=1}^4 \E\left(\sum_{k=t\Delta_n^{-1}+1}^{l_n-t\Delta_n^{-1}-1}Z_{n,k}^{(i)}\right)^{4}
                \leq\mbox{const.}\sum_{i=1}^4\left(l_n\E((Z_{n,1}^{(i)})^4)+l_n^2(\E((Z_{n,1}^{(i)})^2))^2\right)\leq \mbox{const.}(l_n\Delta_n)^2.
    \end{eqnarray*}
    In total we receive
    \begin{eqnarray} \label{A.55}
        \limsup_{n\to\infty}\frac{1}{l_n\Delta_n}\E(|Z_{n,1}^*|^2\1_{\{|Z_{n,1}^*|>\epsilon \sqrt{n\Delta_n}\}})
        \leq  \limsup_{n\to\infty}\frac{1}{l_n\Delta_n}\frac{1}{\epsilon^2n\Delta_n}\E|Z_{n,1}^*|^4=0.
    \end{eqnarray}
  Finally, \eqref{A.113}-\eqref{A.55} result in  the Lindeberg condition.
    \hfill$\Box$\\

    Thus, by \cref{Cov S} and \cref{Lemma:Lindeberg} the assumptions of the central limit theorem of Lindeberg-Feller (see \cite{jacodB}, Theorem VI.5.5.2) are satisfied and we can conclude the
    weak convergence
    \begin{eqnarray} \label{A.CLT}
        \frac{1}{\sqrt{n\Delta_n}}S_{1,n}^*
        \weak \mathcal{N}\left(0,\Sigma(c_1,c_2,c_3)\right) \quad \mbox{ as }n\to\infty.
        \end{eqnarray}
    Moreover, since $(Z_{n,k}^{(i)})_{k\in\N}$ as given in \eqref{def Z} and $(Z_{n,k}^{(5)})_{k\in\N}:=(c_3^T\vect(\xi_{n,k} \xi_{n,k}^T-\Sigma_n))_{k\in\N}$ are uncorrelated
     sequences, $\E((Z_{n,k}^{(i)})^2)<\mbox{const. }\cdot\Delta_n$ and $l_n\Delta_n\to\infty$, we obtain with Markov's inequality
    \begin{eqnarray*}
        \mathbb{P}\left(\frac{1}{\sqrt{n\Delta_n}}|S_{2,n}^*|>\epsilon\right)
        &\leq& \frac{\mbox{const.}}{\epsilon^2}\frac{1}{n\Delta_n}\sum_{i=1}^5 \E\left|\left[\sum_{k=1}^{t\Delta_n^{-1}}
                    +\sum_{i=1}^{\lfloor n/l_n\rfloor-1}\sum_{k=il_n-t\Delta_n^{-1}}^{il_n+t\Delta_n^{-1}}
                +\sum_{k=\lfloor n/l_n\rfloor l_n-t\Delta_n^{-1}}^{n}\right]Z_{n,k}^{(i)}\right|^2\\
        &\leq&\mbox{const. }\frac{1}{n\Delta_n}\frac{n}{l_n\Delta_n}\Delta_n=\frac{1}{l_n\Delta_n}\to 0,
    \end{eqnarray*}
    such that $\frac{1}{\sqrt{n\Delta_n}}S_{2,n}^*\stackrel{\mathbb{P}}{\to}0$. This in combination
    with \eqref{sn*2} and \eqref{A.CLT} result in
    \begin{eqnarray*}
        \frac{1}{\sqrt{n\Delta_n}}S_{n}^*\weak \mathcal{N}\left(0,\Sigma(c_1,c_2,c_3)\right) \quad \mbox{ as }n\to\infty,
    \end{eqnarray*}
    and, in particular, the convergence of the two-dimensional distribution. However,
    we have to be sure that the limit distribution is as stated. From \cref{limits_structure} we already know that \linebreak $sc_1^T\Sigma^*c_1+(t-s)c_2^T\Sigma^*c_2$
    is the covariance matrix of the normally distributed random variable
    \begin{eqnarray*}
        {c_1^T\vect(\Sigma^{1/2}\Wbf_s\Sigma^{1/2 T},\Sigma^{1/2}\Wbf_s^T\Sigma^{1/2 T})
        +c_2^T\vect(\Sigma^{1/2}(\Wbf_t-\Wbf_s)\Sigma^{1/2 T},\Sigma^{1/2}(\Wbf_t-\Wbf_s)^T\Sigma^{1/2 T})}
        =:N^*(c_1,c_2).
    \end{eqnarray*}
     This means
       $ \frac{1}{\sqrt{n\Delta_n}}S_{n}^*\weak N^*(c_1,c_2)+c_3^T\vect(\Wbf^*(\Upsilon)) \mbox{ as }n\to\infty,$
    and the Cram\'{e}r-Wold technique  gives the converges of the two-dimensional distribution
    as stated.
    \hfill$\Box$\\

    To prove the tightness of $\left(\frac{1}{\sqrt{n\Delta_n}}S_n(t)\right)_{t\in[0,T]}$ we use the following criteria
    so that we can apply \cite[Theorem~13.5]{billingsley:1999}.

    \begin{Lemma} \label{Tightness}
    There exists a constant
    $K>0$ such that for any $0\leq r<s<t\leq T$:
    \begin{eqnarray*} \label{eq:1.1}
        \frac{1}{(n\Delta_n)^2}\E(\|S_n(t)-S_n(s)\|^2\|S_n(s)-S_n(r)\|^2)\leq K(t-r)^2.
    \end{eqnarray*}
    \end{Lemma}
    \paragraph{Proof.}
    Without loss of generality $p=1$ and $d=1$, otherwise prove the statement componentwise.
    Therefore, we define
    \begin{eqnarray*}
        V_{n,k}^{(1)}(u_1,u_2):=\sum_{j=\lfloor u_1/\Delta_n\rfloor +1}^{\lfloor u_2/\Delta_n \rfloor}\xi_{n,k-j}
        \quad \mbox{ and } \quad V_{n,k}^{(2)}(u_1,u_2):=\sum_{j=\lfloor u_1/\Delta_n\rfloor +1}^{\lfloor u_2/\Delta_n \rfloor}\xi_{n,k+j} \quad\mbox{ for }0\leq u_1<u_2<\infty,
    \end{eqnarray*}
    such that
    \begin{eqnarray*}
        S_n^{(i)}(t)=\sum_{k=1}^n\xi_{n,k}V_{n,k}^{(i)}(0,t) .
    \end{eqnarray*}
    Note $\E(V_{n,k}^{(i)}(u_1,u_2)^2)\leq \mbox{const.\,} (u_2-u_1)$ and $\E(V_{n,k}^{(i)}(u_1,u_2)^4)\leq \mbox{const.\,} (u_2-u_1)$.
    Moreover,
    \begin{eqnarray} \label{D5}
        \lefteqn{\E((S_n^{(i)}(t)-S_n^{(i)}(s))^2(S_n^{(i)}(s)-S_n^{(i)}(r))^2)}\\
            &&=\sum_{k=1}^n\E\left(\xi_{n,k}^4V_{n,k}^{(i)}(s,t)^2
                V_{n,k}^{(i)}(r,s)^2\right)
                + \sum_{k_1=1}^n\sum_{k_2=1\atop k_2\not=k_1}^{n}    \E\left(\xi_{n,k_1}^2V_{n,k_1}^{(i)}(s,t)^2\xi_{n,k_2}^2V_{n,k_2}^{(i)}(r,s)^2\right). \nonumber
    \end{eqnarray}
We investigate the different summands. First,
\begin{eqnarray} \label{D3}
    \E\left(\xi_{n,k}^4V_{n,k}^{(i)}(s,t)^2
                V_{n,k}^{(i)}(r,s)^2\right)
     =\E\left(\xi_{n,k}^4\right)\E(V_{n,k}^{(i)}(s,t)^2)
                \E(V_{n,k}^{(i)}(r,s)^2)
     \leq \mbox{const.\,}(t-s)(s-r)\Delta_n.
\end{eqnarray}
Next for $|k_1-k_2|>(t-r)/\Delta_n$,
\begin{eqnarray}
    \E(\xi_{n,k_1}^2V_{n,k_1}^{(i)}(s,t)^2
                \xi_{n,k_2}^2V_{n,k_2}^{(i)}(r,s)^2)
                \leq \mbox{const. }(t-r)^2\Delta_n^2,
\end{eqnarray}
and $|k_1-k_2|\leq (t-r)/\Delta_n$, $k_1\not=k_2$,
\begin{eqnarray} \label{D1}
    \E\left(\xi_{n,k_1}^2V_{n,k_1}^{(i)}(s,t)^2
                \xi_{n,k_2}^2V_{n,k_2}^{(i)}(r,s)^2\right)\leq \mbox{const. }[(t-r)+(t-r)^2]\Delta_n^2.
\end{eqnarray}
A conclusion of \eqref{D5}-\eqref{D1} is that
\begin{eqnarray} \label{eq:1.2}
    \E((S_n^{(i)}(t)-S_n^{(i)}(s))^2(S_n^{(i)}(s)-S_n^{(i)}(r))^2)\leq \mbox{const.}(n\Delta_n)^2(t-r)^2.
\end{eqnarray}
Finding an upper bound for $\E((S_n^{(1)}(t)-S_n^{(1)}(s))^2(S_n^{(2)}(s)-S_n^{(2)}(r))^2)$
is alike.
Similar but more technical and tedious calculations as above  yield
\begin{eqnarray} \label{eq:1.3}
    \E((S_n^{(1)}(t)-S_n^{(1)}(s))^2(S_n^{(2)}(s)-S_n^{(2)}(r))^2)\leq \mbox{const. }(n\Delta_n)^2(t-r)^2.
\end{eqnarray}
Then \eqref{eq:1.2} and \eqref{eq:1.3} result in the statement.
\hfill$\Box$\\

\subsubsection{Proofs of the results in \Cref{limit continuous time}}

\paragraph{Proof of \cref{CLT_2}.}
In     \cref{Convergence:fidis} we have already proved the convergence of the two-dimensional distribution.
    The convergence of
    the finite-dimensional distribution is an obvious extension.
     Hence, the weak convergence is a consequence of  \cref{Tightness} and \cite{jacodB}, Theorem VI.4.1 (cf.~\cite{billingsley:1999}, Theorem~13.5).
\hfill$\Box$

\paragraph{Proof of \cref{Lemma1.3}.}
 Let $T>0$. Define $S_n^{(i)}(t):=0_{pd\times pd}$ and $\gbf_l(t):=0_{pd\times pd}$ for $t\leq 0$, $i=1,2$, and $S^{(1)}(t):=\Sigma^{1/2}\Wbf_t\Sigma^{1/2\,T}$
 and $S^{(2)}(t):=\Sigma^{1/2}\Wbf_t^T\Sigma^{1/2\,T}$ for $t\geq 0$. Moreover, for $i=1,2$, $l=1,\ldots,M$, we have
\begin{eqnarray*}
    \int_0^{T}\gbf_l(s)\, S_n^{(i)}(\mbox{d}s)
        =-\left[\int_0^TS_n^{(i)T}(s)\,\mbox{d}\gbf_l(s)^T\right]^T+\gbf_l(T)S_n^{(i)}(T).
\end{eqnarray*}
Applying \cref{CLT_2}, \cite[Theorem VI.6.22]{jacodB}
and partial integration (cf. \cite{Protter:2005}, Corollary II.6.2) we obtain as $n\to\infty$,
\begin{eqnarray*}
    \lefteqn{\left(\frac{1}{\sqrt{n\Delta_n}} \int_0^{T}\gbf_l(s)^T\, S_n^{(i)}(\mbox{d}s)^T,\frac{1}{\sqrt{n\Delta_n}} \int_0^{T}\gbf_l(s)\, S_n^{(i)}(\mbox{d}s)\right)_{l=1,\ldots,M,i=1,2}}\\
    &&\stackrel{\mathcal{D}}{\to}\left(
    -\left[\int_0^TS^{(i)}(s)\,\mbox{d}\gbf_l(s)\right]^T+\gbf_l(T)^TS^{(i)}(T)^T,
    -\left[\int_0^TS^{(i)}(s)^T\,\mbox{d}\gbf_l(s)^T\right]^T+\gbf_l(T)S^{(i)}(T)\right)_{l=1,\ldots,M,i=1,2}\\
        &&=\left(\int_0^T\gbf_l(s)^TS^{(i)}(\mbox{\rm d}s)^T,\int_0^T\gbf_l(s)S^{(i)}(\mbox{\rm d}s)\right)_{l=1,\ldots,M, i=1,2}.
\end{eqnarray*}
Note that
 $\int_0^T \gbf_l(s)S^{(i)}(\mbox{\rm d}s)\stackrel{\p}{\to}\int_0^\infty\gbf_l(s)S^{(i)}(\mbox{\rm d}s)$
and $\int_0^T \gbf_l(s)^TS^{(i)}(\mbox{\rm d}s)^T\stackrel{\p}{\to}\int_0^\infty\gbf_l(s)^TS^{(i)}(\mbox{\rm d}s)^T$
 as $T\to\infty$.
Moreover, for $\epsilon>0$  a conclusion of Markov's inequality is
\begin{eqnarray*}
    \lefteqn{\p\left(\left\|\frac{1}{\sqrt{n\Delta_n}}\sum_{j=\lfloor T/\Delta_n\rfloor+1}^\infty\gbf_l(\Delta_n j)\left[S_n^{(1)}(j\Delta_n)-S_n^{(1)}((j-1)\Delta_n)\right]
        \right\|>\epsilon\right)}\\
        &&\leq \mbox{const. }\frac{1}{n\Delta_n}\sum_{k=1}^n\E\|\xi_{n,k+1}\|^2\left(\sum_{j=\lfloor T/\Delta_n\rfloor+1}^\infty
                \|\gbf_l(\Delta_n j)\|^2\E\|\xi_{n,k-j}\|^2\right)\\
        &&\leq \mbox{const. }\Delta_n\sum_{j=\lfloor T/\Delta_n\rfloor+1}^\infty\|\gbf_l(\Delta_n j)\|^2
        \leq \mbox{const. }\int_T^\infty\|\gbf_l(s)\|^2\,ds \stackrel{T\to\infty}{\to}0.
\end{eqnarray*}
The same statement holds with $S_n^{(2)}$ replaced by $S_n^{(1)}$, and taking the transposed processes.
Hence, the conclusion follows by a convergence together argument (cf. \cite{billingsley:1999}, Theorem~3.2) and the continuous mapping theorem if we take the
transpose of  $\int_0^{\infty}\gbf_l(s)^T\, S_n^{(i)}(\mbox{d}s)^T=[\int_0^{\infty} S_n^{(i)}(\mbox{d}s)\gbf_l(s)]^T$.
\hfill$\Box$

\paragraph{Proof of \cref{Appendix:1}.}
(a) \, Let  $(\gabf_{\Lbf},\Sibf_{\Lbf},\nu_{\Lbf})$ be the characteristic triplet of
$(\Lbf_t)_{t\geq 0}$ and let $\bbs^{m}_\epsilon = \{ \xbf \in \bbr^m : \lVert \xbf \rVert \leq \epsilon \}$ be
a ball around $0_m$ in $\R^m$ with radius $\epsilon>0$. We factorize the L\'{e}vy measure $\nu_{\Lbf}$ into
two L\'{e}vy measures
\beao
\nu_{\Lbf_1^{(\epsilon)}} (A) := \nu_{\Lbf}(A \backslash \bbs^{m}_\epsilon) \quad\quad \mbox{ and
}\quad \quad \nu_{\Lbf_2^{(\epsilon)}} (A): = \nu_{\Lbf}(A \cap \bbs^{m}_\epsilon) \quad \mbox{ for } A\in\mathcal{B}(\R^m\backslash\{0_m\})
\eeao
such that $\nu_{\Lbf}=\nu_{\Lbf_1^{(\epsilon)}}+\nu_{\Lbf_2^{(\epsilon)}}$. Then we can decompose $(\Lbf_t)_{t\geq 0}$ in two
independent L\'{e}vy processes
\beam \label{A1}
    \Lbf_t = \Lbf_1^{(\epsilon)}(t) + \Lbf_{2}^{(\epsilon)}(t) \quad \mbox{ for } t\geq 0,
\eeam
where $\Lbf_i^{(\epsilon)}$ has L\'{e}vy measure $\nu_{\Lbf_i^{(\epsilon)}}$ and expectation $0_m$ ($i=1,2$),
and $\Lbf_1^{(\epsilon)}$ is without Gaussian part.
First, we will show that
\begin{eqnarray} \label{TS}
    \lim_{n\to\infty}\Delta_n^{-1}\E(\Lbf_1^{(\epsilon)}(\Delta_n)\Lbf_1^{(\epsilon)}(\Delta_n)^T\otimes \Lbf_1^{(\epsilon)}(\Delta_n)\Lbf_1^{(\epsilon)}(\Delta_n)^T)
        =\int_{\R^m} \xbf\xbf^T\otimes \xbf\xbf^T\,\nu_{\Lbf_1^{(\epsilon)}}(d\xbf).
\end{eqnarray}
Since the L\'{e}vy measure of $\Lbf^{(\epsilon)}_1$ is finite and $\Lbf^{(\epsilon)}_1$  is without Gaussian part,
$\Lbf^{(\epsilon)}_1$ has the representation as a compound Poisson process with drift
\beam \label{compound}
    \Lbf_1^{(\epsilon)}(t)=\sum_{k=1}^{N(t)}\Jbf_k^{(\epsilon)}+c_1^{(\epsilon)}t, \quad t\geq 0,
\eeam
where $(\Jbf_k^{(\epsilon)})_{k\in\N}$ is a sequence of iid random vectors independent of the Poisson process $(N(t))_{t\geq 0}$
with intensity $\lambda_\epsilon:=\nu_{\Lbf_1^{(\epsilon)}}(\R^m)$.  The distribution
of $\Jbf_k^{(\epsilon)}$ ist $\lambda_\epsilon^{-1}\nu_{\Lbf_1^{(\epsilon)}}$. Moreover, $c_1^{(\epsilon)}$ is a vector
in $\R^{m\times m}$.
We will use on the one hand, that for $l\geq 1$,
\beam   \label{A3a}
    \frac{\p(N(\Delta_n)=l)}{\Delta_n}=\e^{-\lambda_\epsilon \Delta_n}\frac{(\lambda_\epsilon \Delta_n)^l}{\Delta_nl!}\leq \mbox{const. }\p(N(1)=l),
\eeam
and on the other hand, that
\beam \label{A4a}
    \nlim \frac{\p(N(\Delta_n)=l)}{\Delta_n}=\left\{\begin{array}{ll}
    \lambda_\epsilon & \mbox{for } l=1,\\
    0 & \mbox{for } l\geq 2.
    \end{array} \right.
\eeam
Then
\beao
    \lefteqn{\E(\Lbf_1^{(\epsilon)}(\Delta_n)\Lbf_1^{(\epsilon)}(\Delta_n)^T\otimes \Lbf_1^{(\epsilon)}(\Delta_n)\Lbf_1^{(\epsilon)}(\Delta_n)^T)}\\
        &=&\p(N(\Delta_n)=1)\E((\Jbf_1^{(\epsilon)}+c_1^{(\epsilon)}\Delta_n)(\Jbf_1^{(\epsilon)}+c_1^{(\epsilon)}\Delta_n)^T\otimes
            (\Jbf_1^{(\epsilon)}+c_1^{(\epsilon)}\Delta_n)(\Jbf_1^{(\epsilon)}+c_1^{(\epsilon)}\Delta_n)^T)\\
        &&+\sum_{m=2}^\infty\p(N(\Delta_n)=m)\E\left[\left(\sum_{k=1}^m(\Jbf_k^{(\epsilon)}+c_1^{(\epsilon)}\Delta_n)\right)
        \left(\sum_{k=1}^m(\Jbf_k^{(\epsilon)}+c_1^{(\epsilon)}\Delta_n)\right)^T\right.\\
        &&\quad\quad\hspace*{3cm}\left.\otimes\left(\sum_{k=1}^m(\Jbf_k^{(\epsilon)}+c_1^{(\epsilon)}\Delta_n)\right)\left(\sum_{k=1}^m(\Jbf_k^{(\epsilon)}+c_1^{(\epsilon)}\Delta_n)\right)^T\right]\\
        &=:&I_{n,1}+I_{n,2}.
\eeao
Due to \eqref{A3a}, \eqref{A4a} and dominated convergence we get
$
    \lim_{n\to\infty}\Delta_n^{-1}I_{n,2}=0_{m^2\times m^2},
$
and
\begin{eqnarray*}
    \lim_{n\to\infty}\Delta_n^{-1}I_{n,1}=\lambda_\epsilon \E(\Jbf_1^{(\epsilon)}\Jbf_1^{(\epsilon)T}\otimes \Jbf_1^{(\epsilon)}\Jbf_1^{(\epsilon)T})=\int_{\R^m} \xbf\xbf^T\otimes \xbf\xbf^T\,\nu_{\Lbf_1^{(\epsilon)}}(d\xbf),
\end{eqnarray*}
so that \eqref{TS} follows. Moreover, by H\"{o}lder's inequality and \cite[Lemma~3.1]{Asmussen:Rosinski} we have
\begin{eqnarray*}
    &&\hspace*{-0.3cm}\Delta_n^{-1}\|\E(\Lbf_1^{(\epsilon)}(\Delta_n)\Lbf_1^{(\epsilon)}(\Delta_n)^T\otimes \Lbf_1^{(\epsilon)}(\Delta_n)\Lbf_1^{(\epsilon)}(\Delta_n)^T)-
    \E(\Lbf_{\Delta_n}\Lbf_1^{(\epsilon)}(\Delta_n)^T\otimes \Lbf_1^{(\epsilon)}(\Delta_n)\Lbf_1^{(\epsilon)}(\Delta_n)^T)\|\\
    &&\leq [\Delta_n^{-1}\E\|\Lbf_2^{(\epsilon)}(\Delta_n)\|^4]^{1/4}[\Delta_n^{-1}\E\|\Lbf_1^{(\epsilon)}(\Delta_n)\|^4]^{3/4}\\
    &&\stackrel{n\to\infty}{\to}\left[\int_{\R^m} \|\xbf\|^4\,\nu_{\Lbf_2^{(\epsilon)}}(d\xbf)\right]^{1/4}
            \left[\int_{\R^m} \|\xbf\|^4\,\nu_{\Lbf_1^{(\epsilon)}}(d\xbf)\right]^{3/4}
    \stackrel{\epsilon\downarrow 0}{\to}0.
\end{eqnarray*}
On this way we can recursively derive that
\begin{eqnarray*}
    \lefteqn{\Delta_n^{-1}\|\E(\Lbf_1^{(\epsilon)}(\Delta_n)\Lbf_1^{(\epsilon)}(\Delta_n)^T\otimes \Lbf_1^{(\epsilon)}(\Delta_n)\Lbf_1^{(\epsilon)}(\Delta_n)^T)-
    \E(\Lbf_{\Delta_n}\Lbf_{\Delta_n}^T\otimes \Lbf_{\Delta_n}\Lbf_{\Delta_n}^T)\|}\\
    &\leq& \Delta_n^{-1}\|\E(\Lbf_1^{(\epsilon)}(\Delta_n)\Lbf_1^{(\epsilon)}(\Delta_n)^T\otimes \Lbf_1^{(\epsilon)}(\Delta_n)\Lbf_1^{(\epsilon)}(\Delta_n)^T)-
    \E(\Lbf_{\Delta_n}\Lbf_1^{(\epsilon)}(\Delta_n)^T\otimes \Lbf_1^{(\epsilon)}(\Delta_n)\Lbf_1^{(\epsilon)}(\Delta_n)^T)\|\\
    &&\quad +\Delta_n^{-1}\|\E(\Lbf_{\Delta_n}\Lbf_1^{(\epsilon)}(\Delta_n)^T\otimes \Lbf_1^{(\epsilon)}(\Delta_n)\Lbf_1^{(\epsilon)}(\Delta_n)^T)-
    \E(\Lbf_{\Delta_n}\Lbf_{\Delta_n}^T\otimes \Lbf_1^{(\epsilon)}(\Delta_n)\Lbf_1^{(\epsilon)}(\Delta_n)^T)\|\\
    &&\quad +\Delta_n^{-1}\|\E(\Lbf_{\Delta_n}\Lbf_{\Delta_n}^T\otimes \Lbf_1^{(\epsilon)}(\Delta_n)\Lbf_1^{(\epsilon)}(\Delta_n)^T)-
    \E(\Lbf_{\Delta_n}\Lbf_{\Delta_n}^T\otimes \Lbf_{\Delta_n}\Lbf_1^{(\epsilon)}(\Delta_n)^T)\|\\
    &&\quad+\Delta_n^{-1}\|\E(\Lbf_{\Delta_n}\Lbf_{\Delta_n}^T\otimes \Lbf_{\Delta_n}\Lbf_1^{(\epsilon)}(\Delta_n)^T)-
    \E(\Lbf_{\Delta_n}\Lbf_{\Delta_n}^T\otimes \Lbf_{\Delta_n}\Lbf_{\Delta_n}^T)\|\\
    &\stackrel{n\to\infty,\atop\epsilon\downarrow 0}{\to}&0,
\end{eqnarray*}
and hence, the statement follows. \\
(b) \,
When we show that
\begin{eqnarray} \label{A.9}
    \lim_{n\to\infty}\Delta_n^{-1}\E\|\xi_{n,1}\xi_{n,1}^T\otimes \xi_{n,1}\xi_{n,1}^T-(\Bbf\Lbf_{\Delta_n})(\Bbf\Lbf_{\Delta_n})^T\otimes (\Bbf\Lbf_{\Delta_n})(\Bbf\Lbf_{\Delta_n})^T\|=0,
\end{eqnarray}
we can conclude the statement from (a).
Therefore, we use that as $n\to\infty$,
\begin{eqnarray} \label{A.8}
        \Delta_n^{-1}\E\|\xibf_{n,1}-\Bbf\Lbf_{\Delta_n}\|^4\to 0.
\end{eqnarray}
 This we get
from the representation of the components of $\xibf_{n,1}-\Bbf\Lbf_{\Delta_n}$ as
\begin{eqnarray*}
    e_i^T\xibf_{n,1}-e_i^T\Bbf\Lbf_{\Delta_n}=\int_0^{\Delta_n}e_i^T(\e^{-\Lambda s}-I_{pd})\Bbf\,d\Lbf_s
    =\sum_{k=1}^m\int_0^{\Delta_n}e_i^T(\e^{-\Lambda s}-I_{pd})\Bbf e_k\,d L_k(s)\quad (i=1,\ldots,d).
\end{eqnarray*}
Applying  \cite[Lemma 3.2]{Lindner:Cohen:2013} gives
\begin{eqnarray*}
    \lefteqn{\E(e_i^T\xibf_{n,1}-e_i^T\Bbf\Delta\Lbf_{\Delta_n})^4}\\
    &\leq& m^4\sum_{k=1}^m\E\left(\int_0^{\Delta_n}e_i^T(\e^{-\Lambda s}-I_{pd})\Bbf  e_k\, dL_k(s)\right)^4\\
    &\leq&\mbox{const.}\,\sum_{k=1}^m\left[\int_0^{\Delta_n}(e_i^T(\e^{-\Lambda s}-I_{pd})\Bbf  e_k)^4\, ds+\left(\int_0^{\Delta_n}
        (e_i^T(\e^{-\Lambda s}-I_{pd})\Bbf e_k)^2\,ds\right)^2\right]
    =o(\Delta_n),
\end{eqnarray*}
and hence, \eqref{A.8} follows. Using H\"{o}lder's inequality, \eqref{A.8} and $\E\|\Bbf\Lbf_{\Delta_n}\|^4=O(\Delta_n)$ (cf. \cite{Asmussen:Rosinski}, Lemma~3.1)
we obtain
\begin{eqnarray*}
    \lefteqn{\Delta_n^{-1}\E\|\xi_{n,1}(\Bbf\Lbf_{\Delta_n})^T\otimes (\Bbf\Lbf_{\Delta_n})(\Bbf\Lbf_{\Delta_n})^T-(\Bbf\Lbf_{\Delta_n})(\Bbf\Lbf_{\Delta_n})^T\otimes (\Bbf\Lbf_{\Delta_n})(\Bbf\Lbf_{\Delta_n})^T\|}\\
    &&\leq\mbox{const.} \left(\Delta_n^{-1}\E\|\xibf_{n,1}-\Bbf\Lbf_{\Delta_n}\|^4\right)^{1/4}
    (\Delta_n^{-1}\E\|\Bbf\Lbf_{\Delta_n}\|^4)^{3/4}=o(1).
\end{eqnarray*}
Since $\E\|\xi_{n,1}\|^4=O(\Delta_n)$ (cf. ~\cite{Fasen11a}, Proposition~A.1(b)) as well, we obtain recursively
the statement \eqref{A.9}.
\hfill$\Box$

\paragraph{Proof of \cref{Theorem: limit MCARMA}.}
Note that by assumption $\E(\Lbf_1)=0_m$ and $\E\|\Lbf_1\|^4<\infty$. Hence,
on the one hand, \linebreak $\lim_{n\to\infty}\Delta_n^{-1}\E(\xi_{n,k}\xi_{n,k}^T)=\Bbf\Sigma_\Lbf\Bbf^T$ (cf. proof of Proposition~A.1(g) in~\cite{Fasen11a}) and $\E\|\xi_{n,1}\|^4\leq\mbox{const.}\cdot\Delta_n$
(cf. ~\cite{Fasen11a}, Proposition~A.1(b)). In particular, $\E(\xi_{n,k})=0_m$.
Finally, by \cref{Appendix:1}
\begin{eqnarray*}
    \lim_{n\to\infty}\Delta_n^{-1}\E(\xibf_{n,1}\xibf_{n,1}^T\otimes \xibf_{n,1}\xibf_{n,1}^T)
        =\Bbf\otimes\Bbf\left(\int_{\R^m} \xbf\xbf^T\otimes \xbf\xbf^T\,\nu(d\xbf)\right)\Bbf^T\otimes\Bbf^T=\Bbf\otimes\Bbf\cdot\Upsilon\cdot \Bbf^T\otimes\Bbf^T.
\end{eqnarray*}
Then the assumptions of \cref{Lemma1.3} are satisfied and \cref{Lemma1.3} gives the statement.
\hfill$\Box$

\section{Asymptotic behavior of the sample autocovariance function of MCARMA models} \label{Section:ACF CARMA}

In this section we present the main results of this paper starting with the
asymptotic behavior of the sample autocovariance function of a MCARMA process $\Ybf$ as defined in \eqref{CARMA:observation}
driven by the L\'{e}vy process $(\Lbf_t)_{t\in\R}$.
We will assume that $\E\|\Lbf_1\|^2<\infty$ and $\E(\Lbf_1)=0_m$ so that the autocovariance function
$\Gamma_\Ybf(h)=\E(\Ybf_0\Ybf_h^T)$ for $h\in\R$ is well-defined.
The sample autocovariance function is defined as
\begin{eqnarray*}
    \wh\Gamma_n(h)=\frac{1}{n}\sum_{k=1}^{n- h/\Delta_n}(\Ybf_{k\Delta_n}-\ov \Ybf_n)(\Ybf_{k\Delta_n+h}-\ov \Ybf_n)^T \quad\mbox{ for } h\in\{0,\Delta_n,\ldots,(n-1)\Delta_n\},
\end{eqnarray*}
where $\ov \Ybf_n=n^{-1}\sum_{k=1}^{n}\Ybf_{k\Delta_n}$ is the sample mean.
In our first result
we let the sum going to $n$ and neglect the sample mean $\ov \Ybf_n$, i.e., we investigate $\wh\Gamma_n^*(h)
$ as in \eqref{stern gamma}. Afterwards we derive the general
result for the sample autocovariance function.

\begin{Proposition} \label{Theorem:covariance:light}
Let $(\Ybf_t)_{t\in\R}$ be a MCARMA process as defined in \eqref{CARMA:observation} with kernel function
$\fbf$, covariance function $(\Gamma_\Ybf(h))_{h\in\R}$ and driving
L\'{e}vy process $\Lbf$ satisfying $\E\|\Lbf_1\|^4<\infty$,
 $\E(\Lbf_1)=0_m$, $\E(\Lbf_1\Lbf_1^T)=\Sigma_\Lbf=\Sigma_\Lbf^{1/2}\cdot \Sigma_\Lbf^{1/2\,T}$ and
 having L\'{e}vy measure $\nu$.
 Suppose $(\Delta_n)_{n\in\N}$ is a sequence of positive constants
with $\Delta_n\downarrow 0$ and $\lim_{n\to\infty}n\Delta_n=\infty$.
Assume that there exists a sequence of positive constants
$l_n\to\infty$ with  $n/l_n\to\infty$ and $l_n\Delta_n\to\infty$.
Define
\begin{eqnarray*}
    \Upsilon:=\int_{\R^m}\xbf\xbf^T\otimes \xbf\xbf^T\,\nu(d\xbf)\in\R^{m^2\times m^2},
\end{eqnarray*}
and denote by $(\Wbf_t)_{t\geq 0}$  an $\R^{m\times m}$-valued standard Brownian motion independent of
the $\R^{m\times m}$-valued  random matrix $\Wbf^*(\Upsilon)$  with $\vect(\Wbf^*(\Upsilon))\sim\mathcal{N}(0_{m^2},\Upsilon)$.
Let $\mathcal{H}$ be a finite set in $\bigcap_{n\geq n_0}\{0,\Delta_n,2\Delta_n,\ldots,(n-1)\Delta_n\}$, $n_0\in\N$,
and $\wh\Gamma_n^*(h)$
be as in \eqref{stern gamma}. Then as $n\to\infty$,
\begin{eqnarray*}
 \lefteqn{\hspace*{-0.1cm}\left(\sqrt{n\Delta_n}\left(\wh\Gamma_n^*(h)-\Gamma_\Ybf(h)\right)\right)_{h\in\mathcal{H}}
    \stackrel{\mathcal{D}}{\to}\left(\int_0^\infty\fbf(s)\Wbf^*(\Upsilon)\fbf(s+h)^T\,ds\right.}\\
    &&\hspace*{-0.3cm}\left.+\int_0^{\infty}\left[\int_0^{\infty}\fbf(s)\Sigma_\Lbf^{1/2}\,d\Wbf_u\Sigma_\Lbf^{1/2\,T}\fbf(s+u+h)^T\right]
            +\left[\int_0^{\infty}\fbf(s+u-h)\Sigma_\Lbf^{1/2}\,d\Wbf_u^T\Sigma_\Lbf^{1/2\,T}\fbf(s)^T\right]\,ds\right)_{h\in\mathcal{H}}.
\end{eqnarray*}
\end{Proposition}

From this result we get the asymptotic behavior of the sample autocovariance function.

\begin{Theorem} \label{Theorem1:covariance:light}
Let the assumptions of \cref{Theorem:covariance:light} hold.
 Then as $n\to\infty$,
\begin{eqnarray*}
 \lefteqn{\hspace*{0cm}\left(\sqrt{n\Delta_n}\left(\wh\Gamma_n(h)-\Gamma_\Ybf(h)\right)\right)_{h\in\mathcal{H}}
    \stackrel{\mathcal{D}}{\to}\left(\int_0^\infty\fbf(s)\Wbf^*(\Upsilon)\fbf(s+h)^T\,ds\right.}\\
    &&\hspace*{0cm}\left.+\int_0^{\infty}\left[\int_0^{\infty}\fbf(s)\Sigma_\Lbf^{1/2}\,d\Wbf_u\Sigma_\Lbf^{1/2\,T}\fbf(s+u+h)^T\right]
            +\left[\int_0^{\infty}\fbf(s+u-h)\Sigma_\Lbf^{1/2}\,d\Wbf_u^T\Sigma_\Lbf^{1/2\,T}\fbf(s)^T\right]\,ds\right)_{h\in\mathcal{H}}.
\end{eqnarray*}
\end{Theorem}

A consequence is that in the high frequency setting the convergence rate of the sample autocovariance function
is $\sqrt{n\Delta_n}$ which is slower than the classical $\sqrt{n}$ convergence rate
for models in discrete time (cf. \cref{Theorem:ARMA} below).

Moreover note that only for $h\in\bigcap_{n\geq n_0}\{0,\Delta_n,2\Delta_n,\ldots,(n-1)\Delta_n\}$
we received a consistent and asymptotically normally distributed estimator.
If $n\Delta_n^3\to 0$, then $\wh\Gamma_n(\lfloor h/\Delta_n\rfloor \Delta_n)$
is a consistent and asymptotically normally distributed estimator for $\Gamma_\Ybf(h)$ for any $h>0$
as well.

We want to investigate now several special cases where the limit
process has a simpler structure. First, where the driving L\'{e}vy process is an one-dimensional L\'{e}vy process and
second, where the driving L\'{e}vy process of the MCARMA process is a Brownian motion.

\begin{Corollary} \label{Corollary 4.2}
Let the assumptions of \cref{Theorem:covariance:light} hold.
\begin{itemize}
\item[(a)] Let $m=1$ such that $\Lbf=L$, $\Wbf=W$ are one-dimensional processes and let $N$ be a standard normal
random variable independent of $W$. Then as $n\to\infty$,
\begin{eqnarray*}
 \lefteqn{\left(\sqrt{n\Delta_n}\left(\wh\Gamma_n(h)-\Gamma_\Ybf(h)\right)\right)_{h\in\mathcal{H}}}\\
    &&\stackrel{\mathcal{D}}{\to}\left(((\E(L_1^2))^{-1}\E(L_1^4)-3)\Gamma_\Ybf(h)N
    +\int_0^{\infty}\left[\Gamma_\Ybf(u+h)+\Gamma_\Ybf(u-h)\right] \,dW_u\right)_{h\in\mathcal{H}}.
\end{eqnarray*}
\item[(b)] Let $\Lbf$ be a multivariate Brownian motion. Then $\Upsilon=0_{m^2\times m^2}$ and hence, $\Wbf^*(\Upsilon)=0_{m^2\times m^2}$.
\end{itemize}
\end{Corollary}

A different representation of \cref{Theorem1:covariance:light} is by the vector-representation which gives an
explicit description of the limit covariance matrix. However, it is very technical to write it down
for different covariances. For this reason we restrict our attention to a fixed covariance.

\begin{Corollary} \label{Corollary 4.1}
Let the assumptions of \cref{Theorem:covariance:light} hold. Define $$\Sigma_\Ybf(u):=\int_0^\infty\fbf(s+u)\otimes\fbf(s) \,ds
\quad \mbox{ and }\quad \Sigma_\Ybf^*(u):=\int_0^\infty\fbf(s)\otimes\fbf(s+u) \,ds \quad \mbox{ for }u\in\R.$$ Let
$P_{m,m}$ be the Kronecker permutation matrix and $h\geq 0$.
Then $\E(\vect(\Ybf_0\Ybf_h^T))=\E(\Ybf_h\otimes \Ybf_0)=\Sigma_\Ybf(h)\vect(\Sigma_\Lbf)$ and as $n\to\infty$,
\begin{eqnarray*}
    \lefteqn{\hspace*{-0.5cm}\sqrt{n\Delta_n}\left(\vect(\wh\Gamma_n(h)-\Gamma_\Ybf(h))\right)}\\
    &&\stackrel{\mathcal{D}}{\to}
        \mathcal{N}\left(0_{d^2},\Sigma_\Ybf(h)\cdot\Upsilon\cdot\Sigma_\Ybf(h)^T+
        \int_0^\infty\Sigma_\Ybf(u+h)\cdot\Sigma_\Lbf\otimes\Sigma_\Lbf\cdot\Sigma_\Ybf(u+h)^T\,du\right.\\
        &&\hspace*{3cm}+ \int_0^\infty\Sigma_{\Ybf}^*(u-h)\cdot\Sigma_\Lbf\otimes\Sigma_\Lbf\cdot\Sigma_\Ybf^*(u-h)^{T}\,du\\
        &&\hspace*{3cm}+\int_0^\infty \Sigma_\Ybf (u+h)\cdot\Sigma_\Lbf^{1/2}\otimes\Sigma_\Lbf^{1/2}\cdot P_{m,m}
        \cdot\Sigma_\Lbf^{1/2\,T}\otimes\Sigma_\Lbf^{1/2\,T}\cdot\Sigma_\Ybf^{*}(u-h)^T\,du\\
        &&\hspace*{3cm}\left.+ \int_0^\infty\Sigma_{\Ybf}^*(u-h)\cdot\Sigma_\Lbf^{1/2}\otimes\Sigma_\Lbf^{1/2}\cdot P_{m,m}
        \cdot\Sigma_\Lbf^{1/2\,T}\otimes\Sigma_\Lbf^{1/2\,T}\cdot\Sigma_\Ybf(u+h)^{\,T}\,du\right).
\end{eqnarray*}
\end{Corollary}

The advantage of the representation of the limit distribution as in
\cref{Theorem1:covariance:light} is that we are able
to understand the dependence in the model quite well.
For this reason we get several extensions from this including the asymptotic behavior
of cross-covariances
and cross-correlations between the components. The next corollary shows
the behavior of the  cross-covariances for the different components
of a MCARMA process. It is also straightforward to calculate the cross-correlations.

\begin{Corollary}  \label{Corollary 4.4}
Let the assumptions of \cref{Theorem:covariance:light} hold, and denote by $\gamma_i(h)=\E(\Ybf^{(i)}_0\Ybf^{(i)}_h)$
the covariance function of the $i$-th component and by $\gamma_{ij}(h)=\E(\Ybf^{(i)}_0\Ybf^{(j)}_h)$, $h\in\R$, the cross-covariance
function between the $i$-th and the $j$-th component of $(\Ybf_t)_{t\geq 0}$. Furthermore,
\begin{eqnarray*}
    \wh\gamma_n^{(ij)}(h)=e_i^T\wh\Gamma_n(h)e_j=\frac{1}{n}\sum_{k=1}^{n- h/\Delta_n}
        \left(\Ybf_{k\Delta_n}^{(i)}-\ov\Ybf^{(i)}_n\right)\left(\Ybf_{k\Delta_n+h}^{(j)}-\ov\Ybf^{(j)}_n\right) \quad \mbox{ for } h\in\{0,\Delta_n,\ldots,(n-1)\Delta_n\},
\end{eqnarray*}
is the sample cross-covariance
function between the $i$-th and the $j$-th component, and \linebreak $\ov\Ybf^{(i)}_n=e_i^T\ov\Ybf_n=\frac{1}{n}\sum_{k=1}^n\Ybf_{k\Delta_n}^{(i)}$
is the sample mean of the $i$-th  component of $(\Ybf_t)_{t\geq 0}$.
\begin{itemize}
\item[(a)] Then as $n\to\infty$,
\begin{eqnarray*}
    \lefteqn{\hspace*{-0.2cm}\sqrt{n\Delta_n}\left(\wh\gamma_n^{(ij)}(h)-\gamma_{ij}(h)\right)}\\
        &&\hspace*{-0.9cm}\stackrel{\mathcal{D}}{\to}\mathcal{N}\left(0,\int_{\R^m}\left[\int_0^\infty (e_i^T\fbf(s)\xbf)\cdot (e_j^T\fbf(s+h)\xbf) \, ds\right]^2\,\nu(d\xbf)
        +2\int_0^\infty\gamma_{i}(s)\gamma_j(s)+\gamma_{ij}(s+h)\gamma_{ji}(s-h)\,ds\right).
\end{eqnarray*}
\item[(b)] Assume $\Sigma_\Lbf^{-1/2}\Lbf$ has independent and identically distributed components, identically distributed as
$\wt L$. Then as $n\to\infty$,
\begin{eqnarray*}
    \lefteqn{\sqrt{n\Delta_n}\left(\wh\gamma_n^{(ij)}(h)-\gamma_{ij}(h)\right)}\\
        &&\stackrel{\mathcal{D}}{\to}\mathcal{N}\left(0,[\E(\wt L_1^4)-3\E(\wt L_1^2)]\sum_{l=1}^m
        \left[\int_0^\infty(e_i^T\fbf(s)\Sigma_\Lbf^{1/2}e_l)(e_j^T\fbf(s+h)\Sigma_\Lbf^{1/2}e_l)\,ds\right]^2\right.\\
        &&\left.\hspace*{4.5cm}  +2\int_0^\infty\gamma_{i}(s)\gamma_j(s)+\gamma_{ij}(s+h)\gamma_{ji}(s-h)\,ds\right).
\end{eqnarray*}
\end{itemize}
\end{Corollary}

Finally, we want to present Bartlett's formula for a CARMA process.

\begin{Corollary} \label{Bartlett}
Let $(Y_t)_{t\in\R}$ be a CARMA process satisfying the
assumptions of \cref{Theorem:covariance:light}.
\begin{itemize}
\item[(a)] The autocovariance function of $(Y_t)_{t\in\R}$
is denoted by $(\gamma(h))_{h\in\R}$ and the sample autocovariance function by $(\wh\gamma_n(h))_{h\in\{0,\Delta_n,\ldots,(n-1)\Delta_n\}}$.
Then as $n\to\infty$,
\begin{eqnarray*}
    \left(\sqrt{n\Delta_n }\left(\wh\gamma_n(h)-\gamma(h)\right)\right)_{h\in\mathcal{H}}\stackrel{\mathcal{D}}{\to}
    \mathcal{N}\left(0,(m_{s,t})_{s,t\in\mathcal{H}}\right),
\end{eqnarray*}
where
 $   m_{s,t}=((\E(L_1^2))^{-2}\E(L_1^4)-3)\gamma(s)\gamma(t)+\int_{-\infty}^\infty\gamma(u+s)\gamma(u+t)+\gamma(u+s)\gamma(u-t)\,du.$
\item[(b)] The autocorrelation function of $(Y_t)_{t\in\R}$
is denoted by $(\rho(h))_{h\in\R}$ and $\wh\rho_n(h)=\wh\gamma_n(h)/\wh\gamma_n(0)$ 
    for  $h\in\{0,\Delta_n,\ldots,(n-1)\Delta_n\}$
denotes the sample autocorrelation function. Then as $n\to\infty$,
\begin{eqnarray*}
    \left(\sqrt{n\Delta_n }(\wh\rho_n(h)-\rho(h))\right)_{h\in\mathcal{H}}
        \stackrel{\mathcal{D}}{\to}\mathcal{N}(0,(v_{s,t})_{s,t\in\mathcal{H}}),
\end{eqnarray*}
where
\begin{eqnarray} \label{rho}
    \lefteqn{\hspace*{-2cm}v_{s,t}=\int_{-\infty}^{\infty}\rho(u+s)\rho(u+t)-\rho(u-s)\rho(u+t)+2\rho(s)\rho(t)\rho(u)^2} \nonumber\\
        &&\hspace*{-1cm}-2\rho(s)\rho(u)\rho(u+t)-2\rho(t)\rho(u)\rho(u+s)\,du.
\end{eqnarray}
\end{itemize}
\end{Corollary}

\begin{Remark} $\mbox{}$
    \cite{Lindner:Cohen:2013}, Theorem~3.5, derived the asymptotic behavior
    of the sample autocovariance function of a CARMA process sampled at an equidistant time-grid with distance $\Delta>0$.
    We want to compare their and our results. They proved that
    as $n\to\infty$,
    \begin{eqnarray*}
        \frac{1}{\sqrt{n}}\sum_{k=1}^n(Y_{k\Delta}^2-\gamma(0)^2)\weak\mathcal{N}\left(0,[(\E(L_1^2))^{-2}\E(L_1^4)-3](\E(L_1^2))^2\int_0^\Delta f_\Delta(u)^2\,du+2\gamma(0)^2+4\sum_{k=1}^{\infty}\gamma(k\Delta)^2\right),
    \end{eqnarray*}
    where $f_\Delta:[0,\Delta]\to\R$ is defined as $f_\Delta(u)=\sum_{k=-\infty}^\infty f(u+k\Delta)^2$. If we multiply the
    variance of the Gaussian limit distribution by $\Delta$ and let $\Delta\to 0$, then
    \begin{eqnarray}    \label{eq4.1}
        \lefteqn{\hspace*{-3.5cm}\Delta\left[[(\E (L_1^2))^{-2}\E (L_1^4)-3](\E (L_1^2))^2\int_0^\Delta f_\Delta(u)^2\,du+2\gamma(0)^2+4\sum_{k=1}^{\infty}\gamma(k\Delta)^2\right]} \nonumber\\
        &&\hspace*{-3.5cm}\stackrel{\Delta\downarrow 0}{\to} [(\E (L_1^2))^{-2}\E (L_1^4)-3]\gamma(0)^2+4\int_0^\infty\gamma(s)^2\,ds.
    \end{eqnarray}
    Note that the second term $\Delta \cdot2\gamma(0)^2$ converges to $0$.
    The limit result \eqref{eq4.1} is in line with \cref{Bartlett}, since as $n\to\infty$,
    \begin{eqnarray*}
    \frac{\sqrt{\Delta_n}}{\sqrt{n}}\sum_{k=1}^n(Y_{k\Delta_n}^2-\gamma(0)^2)\weak\mathcal{N}\left(0,[(\E (L_1^2))^{-2}\E (L_1^4)-3]\gamma(0)^2+4\int_0^\infty\gamma(s)^2\,ds\right).
    \end{eqnarray*}
    Hence, for the limit distribution it does not matter if first $n\to\infty$ and then $\Delta\to 0$, or $\Delta_n\to0$
    and $n\to\infty$ at the same time. The analogous phenomenon holds also for the sample autocorrelation function
    in the high frequency and the discrete-time setting as given in \cite[Theorem~3.5]{Lindner:Cohen:2013}
    as well.
    \hfill$\Box$
\end{Remark}

\subsection{Proofs}
\paragraph{Proof of \cref{Theorem:covariance:light}.}
 Let us define $\fbf_\Zbf(s):=\e^{-\Labf s}\Bbf\1_{\left[0,\infty\right)}(s)$
and $\Gamma_\Zbf(h)=\E(\Zbf_0\Zbf_h^T)$ for $h\in\R$ with $\Zbf$ as given in \eqref{CARMA:state}.
By the state space representation \eqref{CARMA:observation} we have
the equalities $\fbf(s)=\Ebf\fbf_\Zbf(s)$, $\Gamma_\Ybf(h)=\Ebf\Gamma_\Zbf(h)\Ebf^T$ and
\begin{eqnarray*}
    \sum_{k=1}^n\left(\Ybf_{k\Delta_n}\Ybf_{k\Delta_n+h}^T-\Gamma_\Ybf(h)\right)=
    \Ebf\left(\sum_{k=1}^n\left[\Zbf_{k\Delta_n}\Zbf_{k\Delta_n+h}^T-\Gamma_\Zbf(h)\right]\right)\Ebf^T.
\end{eqnarray*}
Hence, it is sufficient to investigate the asymptotic behavior of $\frac{\sqrt{\Delta_n}}{\sqrt{n}}\sum_{k=1}^n\left(\Zbf_{k\Delta_n}\Zbf_{k\Delta_n+h}^T-\Gamma_\Zbf(h)\right)$.

The proof has a common ground with the proof of \cite[Theorem~3.6]{Fasen11a}.
A multivariate version of the second order Beveridge-Nelson decomposition presented in~\cite[Equation (28)]{PhillipsSolo92}
gives the representation
\beao
    \Zbf_{k\Delta_n}\Zbf_{k\Delta_n+h}^T&=&\sum_{j=0}^{\infty}\e^{-\Labf \Delta_nj}\xibf_{n,k}\xibf_{n,k}^T\e^{-\Labf^T (\Delta_nj+h)}
        +(\Fbf_{n,k-1}^{(1)}(h)-\Fbf_{n,k}^{(1)}(h))+\sum_{r=1}^{\infty}(\Fbf_{n,k,r}^{(2)}(h)+\Fbf_{n,k,-r}^{(2)}(h))\\
        &&\quad +\sum_{r=1}^{\infty}(\Fbf_{n,k-1,r}^{(3)}(h)+\Fbf_{n,k-1,-r}^{(3)}(h)-\Fbf_{n,k,r}^{(3)}(h)-\Fbf_{n,k,-r}^{(3)}(h)),
\eeao
where
\beao
     \Fbf_{n,k}^{(1)}(h)&=&\sum_{j=0}^{\infty}\sum_{s=j+1}^{\infty}\e^{-\Labf \Delta_ns}\xibf_{n,k-j}\xibf_{n,k-j}^T
        \e^{-\Labf^T(\Delta_ns+h)},\\
    \Fbf_{n,k,r}^{(2)}(h)&=&\sum_{j=\max(0,-r-\lfloor h/\Delta_n\rfloor)}^{\infty}\e^{-\Labf \Delta_nj}\xibf_{n,k}\xibf_{n,k-r}^T
        \e^{-\Labf^T(\Delta_n (j+r)+h)},\\
    \Fbf_{n,k,r}^{(3)}(h)&=&\sum_{j=0}^{\infty}\sum_{s=\max(j+1,-r-\lfloor h/\Delta_n\rfloor)}^{\infty}\e^{-\Labf \Delta_ns}\xibf_{n,k-j}\xibf_{n,k-j-r}^T
        \e^{-\Labf^T(\Delta_n (s+r)+h)}.
\eeao
Then
\beam \label{C1}
    \sum_{k=1}^n\Zbf_{k\Delta_n}\Zbf_{k\Delta_n+h}^T\hspace*{-0.2cm}&=&\hspace*{-0.2cm}\sum_{j=0}^{\infty}\e^{-\Labf \Delta_nj}\left(\sum_{k=1}^n\xibf_{n,k}\xibf_{n,k}^T\right)\e^{-\Labf^T (\Delta_nj+h)}
            +(\Fbf_{n,0}^{(1)}(h)-\Fbf_{n,n}^{(1)}(h))\nonumber\\
        &&\hspace*{-0.3cm}+\sum_{k=1}^n\sum_{r=1}^{\infty}(\Fbf_{n,k,r}^{(2)}(h)+\Fbf_{n,k,-r}^{(2)}(h))
            +\sum_{r=1}^{\infty}(\Fbf_{n,0,r}^{(3)}(h)+\Fbf_{n,0,-r}^{(3)}(h)-\Fbf_{n,n,r}^{(3)}(h)-\Fbf_{n,n,-r}^{(3)}(h))  \nonumber\\
        &\hspace*{-0.2cm} =:&J_{n,1}(h)+J_{n,2}(h)+J_{n,3}(h)+J_{n,4}(h).
\eeam
The proof is divided in several parts. We will show the following:
\begin{itemize}
    \item[(i)] ${\displaystyle\left(\sqrt{n\Delta_n}\left(\frac{1}{n}J_{n,1}(h)-\Gamma_\Zbf(h)\right)\right)_{h\in\mathcal{H}}\weak\left(\int_0^\infty\fbf_\Zbf(s) \Wbf^*(\Upsilon)
    \fbf_\Zbf(s+h)^T\,ds\right)_{h\in\mathcal{H}}}$.
    \item[(ii)] ${\displaystyle\frac{\sqrt{\Delta_n}}{\sqrt{n}}J_{n,2}(h) \stackrel{\p}{\longrightarrow} 0_{pd\times pd}}, \,h\in\mathcal{H}$.
    \item[(iii)] ${\displaystyle\left(\frac{\sqrt{\Delta_n}}{\sqrt{n}}J_{n,3}(h)\right)_{h\in\mathcal{H}}\weak
                \left(\int_0^{\infty}\left[\int_0^{\infty}\fbf_\Zbf(s)\Sigma_\Lbf^{1/2}\,d\Wbf_u\Sigma_\Lbf^{1/2\,T}\fbf_\Zbf(s+u+h)^T\right]\,ds\right.}$\\
                ${\displaystyle \left.\quad\hspace*{3.2cm}
            +\int_0^{\infty}\left[\int_0^{\infty}\fbf_\Zbf(s+u-h)\Sigma_\Lbf^{1/2}\,d\Wbf_u^T\Sigma_\Lbf^{1/2\,T}\fbf_\Zbf(s)^T\right]\,ds\right)_{h\in\mathcal{H}}.}$
    \item[(iv)] ${\displaystyle\frac{\sqrt{\Delta_n}}{\sqrt{n}}J_{n,4}(h) \stackrel{\p}{\longrightarrow} 0_{pd\times pd}}, \,h\in\mathcal{H}$.
\end{itemize}
The proof of (ii) and (iv) follows directly from the proof of \cite[Lemma~5.7]{Fasen11a}.\\
\textsl{Proof of (i).} We will use the equality
$\Gamma_\Zbf(h)=\sum_{j=0}^{\infty}\e^{-\Lambda\Delta_n j}\E(\xibf_{n,1}\xibf_{n,1}^T) \e^{-\Lambda^T(\Delta_n j+h)}$
and
\begin{eqnarray} \label{A.4}
    \sqrt{n\Delta_n}\left(\frac{1}{n}J_{n,1}-\Gamma_\Zbf(h)\right)=\Delta_n \sum_{j=0}^{\infty}\e^{-\Labf \Delta_nj}
         \left(\frac{1}{\sqrt{n\Delta_n}}\sum_{k=1}^n[ \xibf_{n,k}\xibf_{n,k}^T-\E(\xi_{n,1}\xi_{n,1}^T)]\right)\e^{-\Labf^T \Delta_nj}\e^{-\Labf^Th}. 
    \end{eqnarray}
An application of \cref{Theorem: limit MCARMA} yields as $n\to\infty$,
\begin{eqnarray} \label{eq:2.1}
    \frac{1}{\sqrt{n\Delta_n}}\sum_{k=1}^n[ \xibf_{n,k}\xibf_{n,k}^T-\E(\xibf_{n,1}\xibf_{n,1}^T)]\weak\Bbf
        \Wbf^*(\Upsilon)\Bbf^T. \nonumber
\end{eqnarray}
Finally,
we denote by $\gbf_n^{(h)}$ and $\gbf^{(h)}$ maps from $\R^{{pd\times pd}}\to \R^{pd\times pd}$ with
\beam \label{gbf}
    \gbf_n^{(h)}(\Cbf)=\Delta_n\sum_{j=0}^{\infty}\e^{-\Labf \Delta_n j}\Cbf\e^{-\Labf^T \Delta_n j}\e^{-h\Labf^T} \quad\mbox{ and }\quad
    \gbf^{(h)}(\Cbf)=\int_0^{\infty}\e^{-\Labf s}\Cbf\e^{-\Labf^T (s+h)}\,\dd s.
\eeam
Since $\gbf_n^{(h)}$ and $\gbf^{(h)}$ are continuous with $\nlim \gbf_n^{(h)}(\Cbf_n)=\gbf^{(h)}(\Cbf)$ for any sequence
$\Cbf_n,\Cbf\in \R^{pd\times pd}$ with $\nlim \|\Cbf_n-\Cbf\|=0$, we can apply a generalized version of the continuous mapping theorem
(cf.~\cite{whitt2002}, Theorem~3.4.4) to obtain as $n\to\infty$,
\begin{eqnarray*}
\left(\frac{\sqrt{\Delta_n}}{\sqrt{n}}[J_{n,1}(h)-\Gamma_\Zbf(h)]\right)_{h\in\mathcal{H}}&=&\left(g_n^{(h)}\left(\frac{1}{\sqrt{n \Delta_n}}\sum_{k=1}^n[\xibf_{n,k}\xibf_{n,k}^T-\E(\xi_{n,1}\xi_{n,1}^T)]\right)\right)_{h\in\mathcal{H}}\\
    &\weak&
         \left(g^{(h)}\left(\Bbf\Wbf^*(\Upsilon)
         \Bbf^T\right)\right)_{h\in\mathcal{H}}.
\end{eqnarray*}
\textsl{Proof of (iii).} An application of \cref{Theorem: limit MCARMA} and a generalized continuous mapping theorem as above gives
\begin{align*}
    &{\left(\frac{\sqrt{\Delta_n}}{\sqrt{n}}\sum_{k=1}^n\sum_{r=1}^{\infty}(\Fbf_{n,k,r}^{(2)}(h)+\Fbf_{n,k,-r}^{(2)}(h))\right)_{h\in\mathcal{H}}}\\
         &\hspace*{0.6cm}=\left(\Delta_n\sum_{j=0}^{\infty}\e^{-\Labf \Delta_nj}
        \left(\frac{1}{\sqrt{n\Delta_n}}\sum_{k=1}^n\left[\sum_{r=1}^{\infty}\xibf_{n,k}\xibf_{n,k-r}^T\e^{-\Lambda^T(h+\Delta_nr)}+
            \sum_{r=1}^{\lfloor h/\Delta_n\rfloor}\xibf_{n,k}\xibf_{n,k+r}^T\e^{-\Labf^T (h-\Delta_n r)}
            \right.\right.\right.\\
            &\hspace*{0.6cm}\quad\quad\quad\quad\quad\quad\quad\quad\quad\quad\left.\left.\left.+\sum_{r=\lfloor h/\Delta_n\rfloor+1}^\infty\e^{-\Labf (\Delta_n r-h)}\xibf_{n,k}\xibf_{n,k+r}^T\right]\right)\e^{-\Labf^T \Delta_n j}\right)_{h\in\mathcal{H}}\\
        &\hspace*{0.6cm}\weak
        \left(\int_0^{\infty}\e^{-\Lambda s}\left[\int_0^{\infty}\Bbf\Sigma_\Lbf^{1/2}\,d\Wbf_u\Sigma_\Lbf^{1/2\,T}\Bbf^T\e^{-\Lambda^T (h+u)}\right]\e^{-\Lambda^T s}\,ds \right.\\
        & \hspace*{0.6cm}\quad +\int_0^{\infty}\e^{-\Lambda s}\left[\int_0^{h}\Bbf\Sigma_\Lbf^{1/2}\,d\Wbf_u^T\Sigma_\Lbf^{1/2\,T}\Bbf^T\e^{-\Lambda^T (h-u)}\right]\e^{-\Lambda^T s}\,ds\\
        &\hspace*{0.6cm} \left. \quad +\int_0^{\infty}\e^{-\Lambda s}\left[\int_h^{\infty}\e^{-\Lambda (u-h)}\Bbf\Sigma_\Lbf^{1/2}\,d\Wbf_u^T\Sigma_\Lbf^{1/2\,T}\Bbf^T\right]\e^{-\Lambda^T s}\,ds\right)_{h\in\mathcal{H}}\\
        &\hspace*{0.6cm}=\left(\int_0^{\infty}\left[\int_0^{\infty}\fbf_\Zbf(s)\Sigma_\Lbf^{1/2}\,d\Wbf_u\Sigma_\Lbf^{1/2\,T}\fbf_\Zbf(s+h+u)^T
        +\int_0^{\infty}\fbf_\Zbf(s+u-h)\Sigma_\Lbf^{1/2}\,d\Wbf_u^T\Sigma_\Lbf^{1/2\,T}\fbf_\Zbf(s)^T\right]\,ds\right)_{h\in\mathcal{H}},
\end{align*}
since $\fbf_\Zbf(s)=0_{d\times m}$ for $s<0$.
This completes the proof of (iii).
\hfill$\Box$

\paragraph{Proof of \cref{Theorem1:covariance:light}.}
In \cite[Theorem~3.1]{Fasen11a} it is stated that $\sqrt{n\Delta_n}\cdot\ov\Ybf_n\stackrel{\mathcal{D}}{\to}\left(\int_0^\infty\fbf(s)\,ds\right)\cdot\mathcal{N}(0_m,\Sigma_\Lbf)$
as $n\to\infty$.
This means in particular that, $\ov\Ybf_n\stackrel{\p}{\to}0_d$ as $n\to\infty$.
Having this  in mind, the rest of the proof goes as in \cite[Proposition 7.3.4]{Brockwelletal1991} for
MA processes in discrete time by proving $$\sqrt{n\Delta_n}\left(\wh\Gamma_n(h)-\wh\Gamma_n^*(h)\right)
\stackrel{\p}{\to}0_{d\times d},$$
and applying our \cref{Theorem:covariance:light}.
\hfill$\Box$

\paragraph{Proof of \cref{Corollary 4.1}.}
In the multivariate case we get the alternative representation of the limit distribution as
\begin{eqnarray*}
    \lefteqn{\hspace*{-4cm}\vect\left(
    \int_0^{\infty}\int_0^{\infty}\fbf(s)\Sigma_\Lbf^{1/2}\,d\Wbf_u\Sigma_\Lbf^{1/2\,T}\fbf(s+u+h)^T\,ds+
    \int_0^{\infty}\int_0^{\infty}\fbf(s+u-h)\Sigma_\Lbf^{1/2}\,d\Wbf_u^T\Sigma_\Lbf^{1/2\,T}\fbf(s)^T\,ds\right)}\\
    &&\hspace*{-3.5cm}=\int_0^{\infty}\int_0^{\infty}\left[\fbf(s+u+h)\otimes\fbf(s) \,ds\right]\cdot\Sigma_\Lbf^{1/2}\otimes \Sigma_\Lbf^{1/2}\,d\vect(\Wbf_u)\\
    &&\hspace*{-3.5cm}\quad+\int_0^{\infty}\int_0^{\infty}\left[\fbf(s)\otimes\fbf(s+u-h) \,ds\right]\cdot\Sigma_\Lbf^{1/2}\otimes \Sigma_\Lbf^{1/2}\,d\vect(\Wbf_u^T).
\end{eqnarray*}
The final statement follows then with  \cref{limits_structure} and \cref{Theorem1:covariance:light}.
\hfill$\Box$

\paragraph{Proof of \cref{Corollary 4.4}.}
(a) \, For the proof we use \cref{Corollary 4.1} and denote by $\fbf_i=e_i^T\fbf$ the $i$-th row of $\fbf$.
The rules for Kronecker products in \eqref{section:matrix calculation} give
\begin{eqnarray*}
    \lefteqn{e_j^T\otimes e_i^T\cdot\left[\int_0^\infty\left[\Sigma_\Ybf^*(u-h)\cdot\Sigma_\Lbf\otimes\Sigma_\Lbf\cdot\Sigma_\Ybf^*(u-h)^T\right]\,du
    \right]\cdot e_j\otimes e_i}\\
    &&=\int_0^\infty\int_0^\infty\int_0^\infty \fbf_j(s)\Sigma_\Lbf \fbf_j(t)^T \cdot \fbf_i(u-h+s)\Sigma_\Lbf \fbf_i(u-h+t)^T\,ds\,dt\,du\\
    &&=2\int_0^\infty\left[\int_0^\infty \left[\int_{t+h}^\infty \fbf_i(u+s)\Sigma_\Lbf \fbf_i(u)^T\,du\right]\fbf_j(t+s)\Sigma_\Lbf \fbf_j(t)^T\,dt\right]\,ds,
\end{eqnarray*}
and on the same way,
\begin{eqnarray*}
    \lefteqn{e_j^T\otimes e_i^T\cdot\left[\int_0^\infty\left[\Sigma_\Ybf(u+h)\cdot\Sigma_\Lbf \otimes\Sigma_\Lbf\cdot\Sigma_\Ybf(u+h)^T\right]\,du
    \right]\cdot e_j\otimes e_i}\\
    &&=2\int_0^\infty\left[\int_0^\infty \left[\int_0^{u+h} \fbf_i(t+s)\Sigma_\Lbf \fbf_i(t)^T\,dt\right] \fbf_j(u+s)\Sigma_\Lbf \fbf_j(u)^T\,du\right]\,ds.
\end{eqnarray*}
Putting both equalities together yields
\begin{eqnarray} \label{C2}
\lefteqn{\hspace*{-1.5cm}e_j^T\otimes e_i^T\cdot\left[\int_0^\infty\left[\Sigma_\Ybf^*(u-h)\cdot\Sigma_\Lbf \otimes\Sigma_\Lbf\cdot\Sigma_\Ybf^*(u-h)^T
    +\Sigma_\Ybf(u+h)\cdot\Sigma_\Lbf\otimes\Sigma_\Lbf\cdot\Sigma_\Ybf(u+h)^T\right]\,du
    \right]\cdot e_j\otimes e_i} \nonumber \\
    &&\hspace*{-1.5cm}=2\int_0^\infty\left[\int_0^\infty \fbf_j(u+s)\Sigma_\Lbf \fbf_j(u)^T\,du\right]\left[\int_0^\infty \fbf_i(t+s)\Sigma_\Lbf \fbf_i(t)^T\,dt\right]\,ds
    =2\int_0^\infty\gamma_j(s)\gamma_i(s)\,ds.
\end{eqnarray}
The equality in \eqref{Permutation_matrix} and analogous calculations as above give
\begin{eqnarray} \label{C3}
    \lefteqn{e_j^T\otimes e_i^T\cdot\left[ \int_0^\infty\left[\Sigma_\Ybf(u+h)\cdot\Sigma_\Lbf^{1/2}\otimes\Sigma_\Lbf^{1/2}\cdot P_{m,m}
        \cdot\Sigma_\Lbf^{1/2\,T}\otimes\Sigma_\Lbf^{1/2\,T}\cdot\Sigma^*_\Ybf(u-h)^{T}\right.\right.} \nonumber\\
        &&\left.\left.\quad+\Sigma_\Ybf^*(u-h)\cdot\Sigma_\Lbf^{1/2}\otimes\Sigma_\Lbf^{1/2}\cdot P_{m,m}
        \cdot\Sigma_\Lbf^{1/2\,T}\otimes\Sigma_\Lbf^{1/2\,T}\cdot\Sigma_\Ybf(u+h)^{T}\right]\,du\right]\cdot e_j\otimes e_i \nonumber \\
    &&=2\int_0^\infty\gamma_{ij}(s+h)\gamma_{ji}(s-h)\,ds
\end{eqnarray}
as well.
Then (a) follows from \cref{Corollary 4.1}, \eqref{C2} and \eqref{C3}.\\
(b) \, is a conclusion of (a) and \cref{Remark_1}(b).
\hfill$\Box$

\paragraph{Proof of \cref{Bartlett}.}
(a) \, can be calculated similarly to \cref{Corollary 4.2}(a). \\
(b) \, The proof can be done step by step as for MA processes in
\cite[Theorem~7.2.1]{Brockwelletal1991} using (a).
\hfill$\Box$

\section{Asymptotic behavior of the sample autocovariance function of MA models} \label{Section:ACF MA}

In Section~\ref{Section:ACF CARMA} we derived the asymptotic behavior of the sample autocovariance
function of a MCARMA process. On a similar way we derive the analogous results for the
sample autocovariance function of a multivariate MA process in discrete time.
The proofs are only slightly different, and
are therefore omitted.
 The first authors who
investigated the asymptotic behavior of the sample autocovariance function for multivariate MA
processes in a very general setup are \cite{Su:Lund}. A difference between their study and our study is
that they define the covariance of two random matrices $\Ubf,\Vbf$ with $\E(\Ubf)=\E(\Vbf)=0_{m\times m}$ as $\Cov_{SL}(\Ubf,\Vbf):=\E(\Ubf\otimes\Vbf)$
where we use $\Cov_{F}(\Ubf,\Vbf):=\E(\vect(\Ubf)\vect(\Vbf)^T)$. The covariance of Su and Lund $\Cov_{SL}(\Ubf,\Ubf)^T=\E(\Ubf^T\otimes\Ubf^T)$
is not necessarily symmetric if $\Ubf$ is not a symmetric random matrix in contrast $\Cov_{F}(\Ubf,\Ubf)^T=\Cov_{F}(\Ubf,\Ubf)$.

A multivariate MA process has the representation
  \begin{eqnarray} \label{multivariate MA}
    \Ybf_k=\sum_{j=0}^\infty \Cbf_{k-j}\xibf_j  \quad \mbox{ for } k\in\Z,
  \end{eqnarray}
   where $(\xi_k)_{k\in\Z}$ is a sequence of iid random vectors in $\R^m$  and $(\Cbf_j)_{j\in\N_0}$ is a sequence
  of deterministic matrices in $\R^{m\times m}$. We will assume that $\E(\xi_1)=0_m$, $\E\|\xi_1\|^2<\infty$ and $\sum_{j=0}^\infty\|\Cbf_j\|^2<\infty$
  so that  the autocovariance function $\Gamma_\Ybf(h)=\E(\Ybf_0\Ybf_h^T)$ for $h\in\Z$ is well-defined.
The sample autocovariance function is defined as
\begin{eqnarray*}
    \wh\Gamma_n(h)=\frac{1}{n}\sum_{k=1}^{n-h}(\Ybf_{k}-\ov \Ybf_n)(\Ybf_{k+h}-\ov \Ybf_n)^T \quad\mbox{ for } h\in\{0,1,\ldots,n-1\},
\end{eqnarray*}
where $\ov \Ybf_n=n^{-1}\sum_{k=1}^{n}\Ybf_{k}$ is the sample mean.
It has the following asymptotic behavior  in analogy to \cref{Theorem1:covariance:light}.

\begin{Theorem} \label{Theorem:ARMA}
Let $(\Ybf_k)_{k\in\Z}$ be a multivariate MA process as defined in \eqref{multivariate MA}
 with noise sequence $(\xibf_k)_{k\in\Z}$ satisfying
   $\E\|\xibf_1\|^4<\infty$, $\E(\xibf_1)=0_m$ and $\Sigma_{\xibf}=\E(\xibf_1\xibf_1^T)$.
   Moreover, we assume that $(\Cbf_j)_{j\in\N_0}$ satisfies
   $\sum_{j=0}^\infty j\|\Cbf_j\|^2<\infty$, $\sum_{j=0}^\infty \|\Cbf_j\|<\infty$ and $\Cbf_j:=0$ for $j<0$.
  Suppose $\Nbf^*(\Upsilon^*)$
    is an $\R^{m\times m}$-dimensional normal random matrix with $\vect(\Nbf^*(\Upsilon^*))\sim\mathcal{N}(0_{m^2},\Upsilon^*)$
    where
    \begin{eqnarray*}
        \Upsilon^*=\E((\xi_1\otimes\xi_1)\cdot (\xi_1\otimes\xi_1)^T)-\E(\xi_1\otimes\xi_1)\E(\xi_1\otimes\xi_1)^T.
    \end{eqnarray*}
    Furthermore, assume $\Nbf^*(\Upsilon^*)$ is independent from the sequence of iid  $\R^{m\times m}$-valued
    random matrices $(\Nbf_r)_{r\in\N}$ with independent standard normally distributed components. Let $\mathcal{H}\subseteq \N_0$
    be a finite set.
    Then as $n\to\infty$,
\begin{eqnarray*}
 \lefteqn{\left(\sqrt{n}\left(\wh\Gamma_n(h)-\Gamma_\Ybf(h)\right)\right)_{h\in\mathcal{H}}}\\
    &&\hspace*{-0.4cm}\stackrel{\mathcal{D}}{\to}\left(\sum_{j=0}^\infty\Cbf_j\Nbf^*(\Upsilon^*)\Cbf_{j+h}^T
    +\sum_{j=0}^{\infty}\left[\sum_{r=1}^{\infty}\Cbf_{j}\Sigma_{\xibf}^{1/2}\Nbf_r\Sigma_{\xibf}^{1/2\,T}\Cbf_{j+r+h}^T\right]
            +\sum_{j=0}^{\infty}\left[\sum_{r=1}^{\infty}\Cbf_{j}\Sigma_{\xibf}^{1/2}\Nbf_r^T\Sigma_{\xibf}^{1/2\,T}\Cbf_{j+h-r}^T\right]\right)_{h\in\mathcal{H}}.
\end{eqnarray*}
\end{Theorem}

\begin{Remark}
The assumption
$\sum_{j=0}^\infty j\|\Cbf_j\|^2<\infty$ is not a necessary assumption. We require this for our way
of proof because they are necessary for the discrete-time versions of $J_{n,1},\ldots,J_{n,4}$ given in \eqref{C1} to be well-defined.
A guess is that $\sum_{j=0}^\infty \|\Cbf_j\|<\infty$ is a sufficient assumption; it is also sufficient in the one-dimensional case.
\hfill$\Box$
\end{Remark}

The vector-representation of this limit result is the following.

\begin{Corollary} \label{corollary arma vec}
Let the assumptions of \cref{Theorem:ARMA} hold.
 Define
 \begin{eqnarray*}
    \Sigma_r:=\sum_{j=0}^{\infty}\Cbf_{j+r}\otimes \Cbf_j \quad\mbox{ and } \quad \Sigma_r^*:=\sum_{j=0}^{\infty}\Cbf_{j}\otimes \Cbf_{j+r} \quad\mbox{ for } r\in\Z.
 \end{eqnarray*}
Let
$P_{m,m}$ be the Kronecker permutation matrix and $h\in\N_0$.
Then $\E(\vect(\Ybf_0\Ybf_h^T))=\E(\Ybf_h\otimes \Ybf_0)=\Sigma_h\vect(\Sigma)$ and as $n\to\infty$,
    \begin{eqnarray*}
        \lefteqn{\hspace*{-0.1cm}\sqrt{n}\vect\left(\wh\Gamma_n(h)-\Gamma_\Ybf(h)\right)}\\
            &&\hspace*{-0.5cm}\weak\mathcal{N}\left(0_{m^2},\Sigma_h\cdot\Upsilon^*\cdot \Sigma_h^T
                +\sum_{r=1}^\infty\left[\Sigma_{r+h}\cdot\Sigma_{\xibf}\otimes\Sigma_{\xibf}\cdot\Sigma_{r+h}^{T}+\Sigma_{r-h}^*
                \cdot\Sigma_{\xibf}\otimes\Sigma_{\xibf}\cdot\Sigma_{r-h}^{*\,T}\right]\right.\\
            &&\hspace*{-1cm}\quad\left.+ \sum_{r=1}^\infty\left[\Sigma_{r+h}\cdot\Sigma_{\xibf}^{1/2}\otimes\Sigma_{\xibf}^{1/2}\cdot
            P_{m,m}\cdot\Sigma_{\xibf}^{1/2\,T}\otimes\Sigma_{\xibf}^{1/2\,T}\cdot\Sigma_{r-h}^{*\,T}+\Sigma_{r-h}^*\cdot\Sigma_{\xibf}^{1/2}\otimes\Sigma_{\xibf}^{1/2}\cdot
            P_{m,m}\cdot\Sigma_{\xibf}^{1/2\,T}\otimes\Sigma_{\xibf}^{1/2\,T}\cdot\Sigma_{r+h}^{T}\right]\right).
   \end{eqnarray*}
 \end{Corollary}

The limit structure in the discrete-time model (\cref{corollary arma vec}) and in the continuous-time model (\cref{Corollary 4.1})
are the same: sums are only replaced by integrals, and $\Upsilon$ by $\Upsilon^*$.

The cross-covariances between the $i$-the and the $j$-th component of $\Ybf$ is presented next.

 \begin{Corollary}
 Let the assumptions of \cref{Theorem:ARMA} hold, and denote by $\gamma_i(h)=\E(\Ybf^{(i)}_0\Ybf^{(i)}_h)$
the autocovariance function of the $i$-th component and by $\gamma_{ij}(h)=\E(\Ybf^{(i)}_0\Ybf^{(j)}_h)$, $h\in\Z$, the cross-covariance
function between the $i$-th and the $j$-th component of $(\Ybf_k)_{k\in\Z}$. Furthermore,
\begin{eqnarray*}
    \wh\gamma_n^{(ij)}(h)=e_i^T\wh\Gamma_n(h)e_j=\frac{1}{n}\sum_{k=1}^{n-h}
        \left(\Ybf_{k}^{(i)}-\ov\Ybf^{(i)}_n\right)\left(\Ybf_{k+h}^{(j)}-\ov\Ybf^{(j)}_n\right) \quad \mbox{ for } h\in\{0,1,\ldots,n-1\},
\end{eqnarray*}
is the sample cross-covariance
function between the $i$-th and the $j$-th component, and  \linebreak $\ov\Ybf^{(i)}_{n}=e_i^T\ov\Ybf_n=\frac{1}{n}\sum_{k=1}^n\Ybf_{k}^{(i)}$
is the sample mean of the $i$-th  component of $(\Ybf_k)_{k\in\Z}$.
 Then as $n\to\infty$,
\begin{eqnarray*}
    \lefteqn{\hspace*{-0.2cm}\sqrt{n}\left(\wh\gamma^{(ij)}_n(h)-\gamma_{ij}(h)\right)}\\
        &&\hspace*{-0.9cm}\stackrel{\mathcal{D}}{\to}\mathcal{N}\left(0,\E\left(\sum_{r=0}^\infty (e_i^T\Cbf_r\xibf)\cdot (e_j^T\Cbf_{r+h}\xibf) \right)^2-3\gamma_{ij}(h)^2
        +\sum_{r=-\infty}^{\infty}\gamma_{i}(r)\gamma_j(r)+\gamma_{ij}(r+h)\gamma_{ji}(r-h)\right).
\end{eqnarray*}
\end{Corollary}

Most results in the literature, with exception of \cite{Su:Lund}, restricted their attention to cross-covariances for either Gaussian
processes or independent processes where the fourth moment part can be neglected. The result presented
here is an extension.

Finally, we present the well-known Bartlett's formula
(see Theorem 7.2.1 and Proposition 7.3.4 in
  \cite{Brockwelletal1991}, and Theorem 6.3.6 and Corollary 6.3.6.1 in \cite{Fuller:book}).
It is a direct consequence of \cref{Theorem:ARMA}.

\begin{Corollary} \label{Corollary 5.2}
Let $(Y_k)_{k\in\Z}$ be an one-dimensional MA process satisfying the assumptions of \cref{Theorem:ARMA}.
\begin{itemize}
   \item[(a)] The autocovariance function of $(Y_k)_{k\in\Z}$ is denoted by $(\gamma(h))_{h\in\Z}$ and $(\wh\gamma_n(h))_{h\in\{0,\ldots,n-1\}}$
   denotes the sample autocovariance function.
   Then as $n\to\infty$,
  \begin{eqnarray*}
    \left(\sqrt{n}(\wh\gamma_n(h)-\gamma(h))\right)_{h\in\mathcal{H}}\weak
        \mathcal{N}\left(0,(m_{s,t})_{s,t\in\mathcal{H}}\right),
  \end{eqnarray*}
  where
 $   m_{s,t}=((\E(\xi_1^2))^{-2}\E(\xi_1^4)-3)\gamma(s)\gamma(t)+\sum_{r=-\infty}^\infty\gamma(r+s)\gamma(r+t)+\gamma(r+s)\gamma(r-t).$
\item[(b)]  The autocorrelation function of $(Y_t)_{t\in\R}$
is denoted by $(\rho(h))_{h\in\Z}$ and $\wh\rho_n(h)=\wh\gamma_n(h)/\wh\gamma_n(0)$
    for  $h\in\{0,\ldots,n-1\}$
denotes the sample autocorrelation function. Then as $n\to\infty$,
\begin{eqnarray*}
    \left(\sqrt{n}(\wh\rho_n(h)-\rho(h))\right)_{h\in\mathcal{H}}
        \stackrel{\mathcal{D}}{\to}\mathcal{N}(0,(v_{s,t})_{s,t\in\mathcal{H}}),
\end{eqnarray*}
where $(v_{s,t})_{s,t\in\mathcal{H}}$ is as in \eqref{rho}.
\end{itemize}
\end{Corollary}

We see that the limit results in the continuous-time model (\cref{Bartlett}(a)) and in the
discrete-time model (\cref{Corollary 5.2}(a)) differ only by changing sums into integrals
and moments of the white noise into moments of the L\'{e}vy process.

\section*{Acknowledgment}
Financial
support by the Deutsche Forschungsgemeinschaft through the research
grant FA 809/2-2 is gratefully acknowledged.

\bigskip

\noindent\textsc{Institute of Stochastics,
Englerstra{\ss}e 2,
D-76131 Karlsruhe, Germany\\ }
\textsc{\emph{Email:} }\href{mailto: vicky.fasen@kit.edu}{vicky.fasen@kit.edu}
{\small\bibliography{Promo}{}}
\bibliographystyle{Econom}

\end{document}